\documentclass[smallextended]{svjour3}

\usepackage{url} 
\usepackage[normalem]{ulem}
\usepackage{blindtext}
\usepackage{lipsum}
\usepackage{amsfonts}
\usepackage{graphicx}
\usepackage{epstopdf}
\usepackage{algorithm}
\usepackage[noend]{algpseudocode}
\usepackage{nccmath}
\usepackage{dsfont}
\usepackage{mathtools}
\usepackage{xcolor}
\DeclareGraphicsExtensions{.eps,.pdf,.png,.jpg}
\usepackage{caption}

\newcommand{\abs}[1]{\left\vert#1\right\vert}
\newcommand{\norm}[1]{\left\lVert#1\right\rVert}
\newcommand{\Prob}{\mathbb{P}}
\newcommand{\X}{\mathbf{X}}

\newcommand{\x}{\mathbf{x}}

\newcommand{\y}{\mathbf{y}}

\renewcommand{\v}{\mathbf{v}}

\newcommand{\z}{\mathbf{z}}

\newcommand{\I}{\mathbf{I}}

\newcommand{\Expectation}{\mathbb{E}}

\newcommand{\convas}{\overset{a.s.}{\longrightarrow}}

\newcommand{\leaveout}[1]{}
\newcommand{\insertme}[1]{#1}

\newcommand{\reals}{\mathbb{R}}

\newcommand{\Ib}{\textbf{I}}
\newcommand{\fb}{\textbf{f}}
\newcommand{\bP}{\textbf{P}}
\renewcommand{\P}{\mathbf{P}}
\newcommand{\Sigmab}{\mathbf{\Sigma}}
\newcommand{\zero}{\mathbf{0}}
\newcommand{\bA}{\textbf{A}}

\usepackage{enumitem}
\spnewtheorem{assumption}{Assumptions}{\bf}{\rm}
\usepackage{bm,psfrag}

\begin{document}

\title{A stochastic subspace approach to gradient-free optimization in high dimensions}

\author{ David Kozak
\and Stephen Becker
\and Alireza Doostan
\and Luis Tenorio}

\institute{David Kozak
\at Solea Energy. Work completed while in Department of Applied Mathematics and Statistics, Colorado School of Mines, Golden, CO
\and
Stephen Becker \at Department of Applied Mathematics, University of Colorado, Boulder, CO
\and 
Alireza Doostan \at Aerospace Engineering Sciences Department, University of Colorado, Boulder, CO
\and Luis Tenorio \at  Department of Applied Mathematics and Statistics, Colorado School of Mines, Golden, CO
}

\maketitle

\begin{abstract}
We present a stochastic descent algorithm for unconstrained optimization that is particularly efficient when the objective function is slow to evaluate and gradients are not easily obtained, as in some PDE-constrained optimization and machine learning problems.  
The algorithm maps the gradient onto a low-dimensional random subspace of dimension $\ell$ at each iteration, similar to coordinate descent but without restricting directional derivatives to be along the axes. 
Without requiring a full gradient, this mapping can be performed by computing $\ell$ directional derivatives (e.g., via forward-mode automatic differentiation).
We give proofs for convergence in expectation under various convexity assumptions as well as probabilistic convergence results under strong-convexity. 
Our method \leaveout{extends the} \insertme{provides a novel extension to the} well-known Gaussian smoothing technique to descent in subspaces of dimension greater than one, opening the doors to new analysis of Gaussian smoothing when more than one directional derivative is used at each iteration.
We also provide a finite-dimensional variant of a special case of the Johnson-Lindenstrauss lemma. 
Experimentally, we show that our method compares favorably to coordinate descent, Gaussian smoothing, gradient descent and BFGS (when gradients are calculated via forward-mode automatic differentiation) on problems from the machine learning and shape optimization literature.

\keywords{Randomized methods \and gradient-free\and Gaussian processes\leaveout{\and shape optimization}\and stochastic gradients}

\subclass{90C06\and 93B40\and 65K10}

\end{abstract}

\section{Introduction}

We consider optimization problems of the form 
\begin{equation}\label{eq: minimize f}
 \min_{\x \in \reals^d} f(\x),
\end{equation}
where $f : \reals^d \to \reals$ has $\lambda$-Lipschitz gradient but $\nabla f(\x)$ is costly to evaluate. 
We also consider additional restrictions on $f$ such as convexity or $\gamma$-strong convexity, which will be made clear as required. 
The main idea is straightforward and has a long history: descend along directions in input space rather than along the gradient.

Directional derivatives can be obtained exactly by forward-mode automatic differentiation, as discussed in \cite{nesterov2017random}, at a cost of approximately one function evaluation per direction. The gradient can be obtained by performing $d$ such calculations in orthogonal directions.
Reverse-mode automatic differentiation would enable calculation of the gradient at a cost of roughly four function evaluations \cite{nesterov2017random} but it has a potential explosion of memory when creating temporary intermediate variables. 
For example, in unsteady fluid flow, the naive adjoint state method requires storing the entire time-dependent PDE-solution~\cite{NASA_NielsenDiskin}. 
Hybrid check-pointing schemes~\cite{wang2009minimal}, designed to reduce memory-overhead, are the subject of active research but the issue has not yet been satisfactorily resolved.
We desire methods that can make progress towards the optima after fewer than $d$ function evaluations per iteration, while still providing convergence guarantees similar to  those of traditional methods.
To this end, we approximate $\nabla f(\x_k)$ with $\ell$ directional derivatives determined by a random matrix $\P_k \in \reals^{d \times \ell}$. 
Such a choice amounts to descending in an $\ell$-dimensional subspace of gradient space and
results in the following recursion,
\begin{equation}\label{eq: iterations}
    \x_{k+1} = \x_k - \alpha \P_k\P_k^\top \nabla f(\x_k),
\end{equation}
where $\alpha >0$ is fixed, $\P_k \in \reals^{d\times \ell}$ is a random matrix with the properties $\Expectation~\P_k\P_k^\top  = \I_d$ and $\P_k^\top\P_k = (d/\ell)\,\I_{\ell}$. 
Note that when $\P_k\P_k^\top $ is diagonal \eqref{eq: iterations} reduces to randomized block-coordinate descent. 
In this document we show that randomized block-coordinate descent is suboptimal for algorithms of the form \eqref{eq: iterations} due to its strong dependence on both the ambient dimension of the problem and the structure of the gradient.
Using a variant of the Johnson-Lindenstrauss lemma we provide non-asymptotic, probabilistic convergence results with spherically symmetric random matrices $\P_k$, results that we show do not hold for coordinate descent.

For concreteness consider the matrix $\P$ comprised of columns $\P^1, \ldots, \P^\ell$. Then an $\ell$-dimensional subspace approximating the gradient can be obtained using finite-differences
\begin{equation}\label{eq: haar}
    \nabla f(\x) \approx
    \P \begin{pmatrix}
    \frac{f(\x + \P^{1}h) -f (\x)}{h} \\ \vdots \\
    \frac{f(\x + \P^{\ell}h) -f (\x)}{h}
    \end{pmatrix}.
\end{equation}
By using exact directional derivatives obtained with forward-mode automatic differentiation, \eqref{eq: haar} reduces to $\nabla f(\x) \approx \P\P^\top \nabla f(\x)$, resulting in the form for \eqref{eq: iterations}. In this paper we analyze the effect that the choice of matrices $\P$ can have on the convergence of \eqref{eq: iterations}. This is accomplished, in part, by analyzing how well $\P\P^\top \nabla f(\x)$ approximates the gradient.

A particular case of \eqref{eq: minimize f} is Empirical Risk Minimization (ERM) commonly used in machine learning, where $ f(\x)=(1/n)\sum_{i=1}^n f_i(\x)$ and $n$ is typically very large.
Hence an ERM problem is amenable to iterative stochastic methods that approximate $\nabla f(\x)$ using $S$ randomly sampled observations, $(i_s)_{s=1}^{S} \subset \{1, \ldots, n\}$, at each iteration with $f_S(\x) = (1/S)\sum_{s=1}^S f_{i_s}(\x)$ where $S\ll n$. 
While the methods we discuss do not require a finite-sum structure, they can be used for such problems.

There are important classes of functions that do not fit into the ERM framework and therefore do not benefit from stochastic gradient descent which is tailored to ERM. 
Partial Differential Equation (PDE) constrained optimization is one such example, and except in special circumstances (such as \cite{haber2012effective}), a stochastic approach leveraging the ERM structure (such as stochastic gradient descent and its variants) does not provide any benefits.
This is because in PDE-constrained optimization the cost of evaluating each $\nabla f_i(\x)$ is often identical to the cost of evaluating $\nabla f(\x)$. 
Problems outside of the ERM framework are not limited to parameter estimation for PDEs.  For example, parameter estimation of Gaussian processes, specifically the sparse Gaussian process framework of \cite{snelson2006sparse,titsias2009variational} does not benefit from an ERM structure but can benefit from our methodology.

\paragraph{PDE-constrained optimization}
Partial differential equations are frequently used to model physical phenomena.
Successful application of PDEs to modeling is contingent upon appropriate discretization and parameter estimation.
Parameter estimation in this setting arises in optimal control, or whenever the parameters of the PDE are unknown, as in inverse problems.
Algorithmic and hardware advances for PDE-constrained optimization have allowed for previously impossible modeling capabilities.
Examples include fluid dynamics models with millions of parameters for tracking atmospheric contaminants \cite{Flath2011}, modeling the flow of the Antarctic ice sheet \cite{isaac2015scalable,petra2014computational}, parameter estimation in seismic inversion \cite{Abacioglu2001,bui2013computational}, groundwater hydrology \cite{WRCR:WRCR23173}, experimental design \cite{horesh2010optimal,haber2012numerical}, and atmospheric remote sensing \cite{LikelihoodInformedDimensionReduction_Alessio}. 

\paragraph{Gaussian processes} 
Gaussian processes are an important class of stochastic processes. 
In this paper we use them to model an unknown function in the context of regression. 
The celebrated representer theorem of Kimeldorf and Wahba \cite{kimeldorf1970correspondence} allows the modeling of functions from an infinite-dimensional reproducing kernel Hilbert space using only machinery from finite-dimensional linear algebra. However, the applications of Gaussian processes are somewhat hamstrung in many modern settings because their time complexity scales as $\mathcal{O}(n^3)$ and their storage as $\mathcal{O}(n^2)$. 
One recourse is to approximate the Gaussian process, allowing time complexity to be reduced to $\mathcal{O}(nm^2)$ with storage requirements of $\mathcal{O}(nm)$, where $m \ll n$ is the number of points used in lieu of the full data set.
Methods have been developed to place these $m$ inducing points, also called landmark points, along the domain at points different from the original inputs \cite{snelson2006sparse,titsias2009variational};  optimal placement of the landmarks is a continuous optimization problem with dimension equal to the number of inducing points to be placed in addition to the number of parameters to be estimated. 
Such a framework places a great burden on the optimization procedure as improperly placed landmark points may result in poor approximations.
\subsection{Related work}
Despite being among the easiest to understand and oldest variants of gradient descent, subspace methods (by far the most common of which is coordinate descent) have, until recently, attracted relatively little attention in the optimization literature. 

\paragraph{Coordinate descent schemes}
The simplest variant of subspace descent is a deterministic method that cycles over the coordinates.
This method is popular because many problems have structure that makes a coordinate update very cheap.
However convergence results for coordinate descent require challenging analysis and the class of functions for which it converges is restricted; indeed, \cite{warga1963minimizing,powell1973search} provide simple examples for which the method fails to converge while simpler-to-analyze methods such as gradient descent converge. 

Choosing the coordinates randomly can lead to results on par with gradient descent \cite{nesterov2012efficiency,richtarik2014iteration}. 
Much emphasis has been placed recently on accelerating coordinate descent methods \cite{allen2016even,hanzely2018accelerated}, but the improvements require knowledge of the Lipschitz constants of the partial derivatives of the functions and/or special structure in the function to make updates inexpensive and to choose a sampling scheme.
See \cite{wright2015coordinate} for a survey of recent results.

A generalization of coordinate descent for linear systems is provided by \cite{gower2015stochastic} wherein the goal is to solve the dual problem.
The idea proposed in \cite{gower2015stochastic} of descending in a random direction according to some pre-specified distribution that is not uniform makes it more similar to ours than other algorithms that focus on solving the dual problem such as, e.g., \cite{SDCA}. 

\paragraph{Zeroth-order optimization}\label{sec:DFO}
Our methods use directions $\bP^\top\nabla f(\x)$, where $\bP$ is $d\times \ell$ with $\ell \ll d$, which is equivalent to taking $\ell$ directional derivatives of $f$ at $\x$ \insertme{. Observe that as our methods do not use gradients they fall into the class of gradient-free optimization, however since we use exact directional derivatives the methods are not derivative-free.} To be clear, when $\nabla f(\x)$ is readily available, zeroth-order optimization methods are not competitive with first- or second-order methods. 
For example, if $f(\x) = \|\bA\x-\textbf{b}\|^2$, \insertme {with $\bA \in \reals^{n \times d}$ and $\textbf{b} \in \reals^n$} then evaluating $f(\x)$ and evaluating $\nabla f(\x) = 2\bA^\top(\bA\x-\textbf{b})$ have nearly the same computational cost, namely $\mathcal{O}(nd)$.  
In fact, such a statement is true regardless of the structure of $f$: by using reverse-mode automatic differentiation (AD), one can theoretically evaluate $\nabla f(\x)$ in about four-times the cost of evaluating $f(\x)$, regardless of the dimension $d$ \cite{griewank2008evaluating}. 
In the context of PDE-constrained optimization, the popular adjoint-state method, which is a form of AD applied to either the continuous or discretized PDE, also evaluates $\nabla f(\x)$ in time independent of the dimension.
However, there are many situations when AD and the adjoint-state method are inefficient or not applicable. 
Finding the adjoint equation requires a careful derivation (which depends on the PDE as well as on initial and boundary conditions), and then a numerical method must be implemented to solve it, which takes considerable development time. 
For this reason complicated codes that are often updated with new features, such as weather models, rarely have the capability to compute a full gradient.
There are software packages that solve for the adjoint automatically, or run AD, but these require a programming environment that restricts the user, and may not be efficient in parallel high-performance computing environments.

 There is a plethora of derivative-free optimization (DFO) algorithms, including 
 grid search, Nelder-Mead, (quasi-) Monte-Carlo sampling,  simulated annealing and MCMC methods~\cite{SimulatedAnnealing}. 
 Modern algorithms include randomized methods, Evolution Strategies (ES) such as CMA-ES~\cite{CMA-ES2001}, 
 Hit-and-Run~\cite{Hit-And-Run} and random cutting planes~\cite{RandomCuttingPlane}.
 Textbook DFO methods (\cite[Algo.\ 10.3]{connDFO}, 
\cite[Algo.\ 9.1]{NocedalWright}) are based on interpolation and trust-regions. 
A limitation of all these methods is that they do not scale well to high-dimensions (beyond $\mathcal{O}(10^2))$.

\paragraph{Stochastic gradient-free methods}\label{sec:SDFO}
Our stochastic subspace descent (SSD) method \eqref{eq: iterations} has been previously explored under the names ``random gradient,'' ``random pursuit,'' ``directional search'', and ``random search''.
The algorithm dates back to the 1970s, 
  with some analysis (cf.\ 
\cite[Ch.\ 6]{ErmolievWets1988} and
\cite{GavianoRandomSearch1975,SolisWets81}), but it never achieved prominence because zeroth-order methods are not competitive with first-order methods when the gradient is available. 
Most analysis has focused on the specific case $\ell=1$ \cite{nesterov2017random,stich2013optimization,LewisDirectionalSearch}.
More recently, the random gradient method has seen renewed interest. For example, \cite{LewisDirectionalSearch} analyzes the case when $f$ is quadratic, and \cite{stich2013optimization} provides an analysis (assuming a line search oracle).
The method of Gaussian smoothing introduced in  \cite{nesterov2017random} is similar to what we propose. We compare the analysis and performance of \cite{nesterov2017random} to that of our method in Sections \ref{sect: MainResults} and \ref{sect: EmpiricalResults}. \leaveout{Gaussian smoothing considers smoothing a function to make it more tractable} \insertme{Gaussian smoothing convolves the objective function with a Gaussian random variable to make the objective differentiable without changing its stationary points.} 
\begin{equation}\label{eq: gaussiansmoothing}
    f^h(\x) = \Expectation_{\mathbf{u}} f(\x+\mathbf{u}h),
\end{equation}
for $h\geq 0$ and $\mathbf{u}\sim\mathcal{N}(\mathbf{0},\mathbf{\Sigma}_d)$. It is common (e.g., \cite{berahas2019theoretical,berahas2019derivative}), and simpler, to consider the case $\mathbf{\Sigma}
_d = \I_d$. It is shown in \cite{nesterov2017random} that \eqref{eq: gaussiansmoothing}  leads to the following finite-difference approximation of the gradient,
\begin{equation}\label{eq: gaussiansmoothing fd}
    \nabla f(\x) \approx \nabla f^h(\x)= \Expectation_{\mathbf{u}}\left[ \mathbf{u} \frac{f(\x + \mathbf{u} h) - f(\x)}{h}\right].
\end{equation}
The obvious way to estimate $\nabla f^h(\x)$ is the single-sample unbiased estimator proposed by Nesterov,
\begin{equation*}
     \nabla f^h(\x) \approx \mathbf{u} \frac{f(\x + \mathbf{u} h) - f(\x)}{h}.
\end{equation*}
Naturally, such an estimator may have a large variance. Thus, to reduce the variance it is tempting to consider taking $\ell >1$ and averaging the results as follows
\begin{equation}\label{eq: berahas gaussian smoothing}
    \frac{1}{\ell}\sum_{i=1}^\ell \mathbf{u}_i \frac{f(\x+\mathbf{u}_ih)-f(\x)}{h},
\end{equation}
where $\mathbf{u}_1, \ldots, \mathbf{u}_\ell \overset{iid}{\sim} \mathcal{N}(\mathbf{0},\I_d)$, as in, e.g.,  \cite{berahas2019theoretical}. While independent directional derivatives provide an estimate of the gradient with a reduced variance compared to Gaussian smoothing, 
independence comes with the undesirable property highlighted in \cite{berahas2019theoretical}: even $\ell>d$ directional derivatives are insufficient to recover the exact gradient. 
 In this paper we consider an alternative to \eqref{eq: berahas gaussian smoothing} for approximation of the gradient when $\ell\geq 1$ and $h \to 0$, which is valid when using a derivative oracle such as forward-mode automatic differentiation. \insertme{Rather than independent Gaussian vectors, we require the $\mathbf{u}_i$ to be orthonormal; equation \eqref{eqn: Haar} provides a method for generating such vectors. This is discussed further in Section \ref{sect: non-asymptotic}.}
\leaveout{In so doing we provide a generalization of Gaussian smoothing as a mapping of the gradient onto a lower dimensional subspace. This interpretation} \insertme{The use of orthonormal $\mathbf{u}_i$} enables the use of machinery that provides sharper and simpler analysis than previously available.
Various proximal, acceleration and noise-tolerant extensions and analyses of Gaussian smoothing have appeared in \cite{dvurechensky2017randomized,dvurechensky2018accelerated,GhadimiLan13,berahas2019theoretical}.
Another variant of random gradient has recently been proposed in the reinforcement learning community. 
The Google Brain Robotics team sampled orthogonal directions to train reinforcement learning systems~\cite{choromanski2018structured} but treated it as a heuristic to approximate Gaussian smoothing 
Similarly,  \cite{duchi2015optimal} uses \insertme{columns from} Haar-distributed matrices and considers the case $\ell=1$, focusing on technical issues related to the small bias introduced by estimation of directional derivatives by finite differences.
The recent papers \cite{cartis2018global} and \cite{berahas2019theoretical} also investigate techniques similar to ours though like \cite{duchi2015optimal} they focus on the implications of the finite difference bias. 
Following \cite{nesterov2012efficiency} we assume that directional derivatives are available via an oracle such as forward-mode automatic differentiation so the finite-difference bias is of no concern. \insertme{Analysis using finite-differences in place of exact directional derivatives is possible (see, e.g., \cite{berahas2019theoretical}). In a forthcoming manuscript we show that in the case $h>0$, convergence of our algorithm is to within a ball of radius $\mathcal{O}(h^2)$ of the minimum function evaluation where $h$ can be on the order of $10^{-8}$; $h$ smaller than $10^{-8}$ incurs numerical instability errors and should be avoided. For the types of problems discussed in this work, in particular PDE-constrained optimization, precision of this order is often impossible; thus, the error attributable to a biased estimate of the gradient is subsumed by other sources of error such as measurement error or termination of the optimization algorithm prior to convergence. For this reason we omit analysis of the finite-difference case and focus only on the setting of exact directional derivatives.}

\paragraph{Alternatives}
As a baseline one could use $\mathcal{O}(d)$ function evaluations to obtain $\nabla f(\x)$ using forward-mode automatic differentiation, which is too costly when $d$ is large and evaluating $f(\x)$ is expensive.
Once $\nabla f(\x)$ is computed, one can run gradient descent, accelerated variants~\cite{Nesterov1983}, non-linear conjugate gradient methods \cite{cg851}, or quasi-Newton methods like BFGS and its limited-memory variant~\cite{NocedalWright}.
In the numerical results section we compare to (finite-difference versions of) gradient descent and BFGS because they are so ubiquitous. We also provide comparisons to Gaussian smoothing and to coordinate descent as the method we propose generalizes both concepts.

 \subsection{Structure of this document and contributions}
In Section \ref{sect: SSD} we investigate convergence of the stochastic subspace descent method for smooth functions. Assumptions used throughout the document are listed, and expected rates of convergence are provided in the case of non-convex, convex, and strongly-convex functions, as well as functions satisfying the Polyak-Lojasiewicz inequality.  
In Section \ref{sect: non-asymptotic} we discuss the properties of gradient approximation along random orthogonal directions for use with \eqref{eq: iterations}. Choosing directions from a specific distribution that we specify, we are able to provide non-asymptotic, high-probability convergence results for strongly-convex functions. \insertme{As previously mentioned this algorithm is a generalization of several classical algorithms for which convergence has already been studied, however as stochastic subspace descent has not been previously introduced, the main results are original. In particular, Theorem 1 and Corollary 1 provide a generalization and different analysis than has previously been performed on Gaussian/spherical smoothing. Theorem 2 is a straightforward generalization of a previously known result. Theorem 3 is a generalization of known analysis for convergence of gradient descent on non-convex objectives. Theorem 4 is new analysis. Lemma 1 is probably known but we have not seen it stated as such, and the remarks are known or simple to prove. }

In Section \ref{subsect:synthetic-data} we provide empirical results on a simulated function that Nesterov dubs ``the worst function in the world'' \cite{nesterov2013introductory}. 
In Section \ref{subsect: experiments-GP} the placement of inducing points for sparse Gaussian processes in the framework of \cite{titsias2009variational} is optimized. 
As a final empirical demonstration, in Section \ref{subsect: experiments-plate}  our algorithms are tested in the PDE-constrained optimization setting on a shape optimization problem.
For the sake of readability, proofs are relegated to the Appendix. 

In this document, uppercase boldfaced letters represent matrices, lowercase boldfaced letters are vectors. 
The vector norm is assumed to be the Euclidean 2-norm, and the matrix norm is the operator norm. 

\section{Main results}\label{sect: MainResults} 
For the remainder of this section we make use of the following assumptions on the sequence of matrices $(\mathbf{P}_k)$ and the function $f$ to be optimized. 
\begin{assumption} Let $\ell \leq d$ and assume:
\begin{enumerate} 
    \item[(A0)\,\,]  $\mathbf{P}_k \in \reals^{d \times \ell}$, $k=1,2, \ldots$, are iid random matrices such that $\Expectation~ \mathbf{P}_k\mathbf{P}_k^\top  = \I_d$  and $\mathbf{P}^\top _k\mathbf{P}_k=(d/\ell)~\I_{\ell}$.
\end{enumerate}
\begin{enumerate}[label=A\theenumi]
	\item[\refstepcounter{enumi}(A\number\value{enumi})\,\,]    $f: \reals^d \to \reals$ is continuously-differentiable with a $\lambda$-Lipschitz first derivative.
	\item[\refstepcounter{enumi}(A\number\value{enumi})\,\,]The function $f$ attains its minimum $f_*$.
	\item[\refstepcounter{enumi}(A\number\value{enumi})\,\,] For some $0<\gamma \leq \lambda$ (where $\lambda$ is the Lipschitz constant in (A1)) and all $\x \in \reals^d$, the function $f$ satisfies the Polyak-Lojasiewicz (PL) inequality: 
			\begin{equation}\label{eq:PL}
			f(\x) - f_* \leq  \norm{\nabla f(\x)}^2/(2\gamma).
			\end{equation} 
    \item[(\theenumi')]   $f$ is $\gamma$-strongly-convex for some $\gamma > 0$ and all $\x \in \reals^d$.
    Note, $\lambda \geq \gamma$ where $\lambda$ is the Lipschitz constant in (A1).
    \item[(\theenumi'')]   $f$ is convex and attains its minimum $f_*$ on a domain $\mathcal{D}$, and there is an $R>0$ such that for the parameter initialization $\x_0$,
    $\max_{\x,~\x_* \in \mathcal{D}} \{\norm{\x-\x_*} : f(\x) \leq f(\x_0) \} \leq R$.
\end{enumerate}
\end{assumption}
\insertme{The assumptions on the matrix $\P_k$ can be satisfied by sampling $\ell \leq d$ columns without replacement from an orthogonal matrix.} Coercivity of $f$ implies the existence of the constant $R$ in (A3''). 
For the results below, particularly the rate in Theorem \ref{thm: convergence-convex}, we require knowledge of the value of $R$.
Also note that (A3') implies (A3).

\subsection{Asymptotic results} \label{sect: SSD}
We now provide conditions under which function evaluations $f(\x_k)$ of stochastic subspace descent converge to a function evaluation at the optimum $f(\x_*)$. 
In the case of a unique optimum we also provide conditions for the iterates $\x_k$ to converge to the optimum $\x_*$. 
Stochastic subspace descent, so-called because at each iteration the method descends in a random low-dimensional subspace, is a gradient-free method as it only requires computation of directional derivatives at each iteration without requiring direct access to the gradient.
In practice we use $\ell$ columns from a scaled Haar-distributed random matrices to define randomized directions along which to descend at each iteration. 
However, neither Theorem \ref{thm:convergence}, nor the subsequent theorems in this subsection require Haar-distributed matrices specifically, as long as the random matrices satisfy Assumption (A0). 
Section \ref{sect: non-asymptotic} demonstrates the advantages of using Haar over random coordinate descent type schemes.

\begin{theorem}[Convergence of SSD]\label{thm:convergence}
 Assume (A0), (A1), (A2), (A3) and let $\x_0$ be an arbitrary initialization. 
 Then recursion \eqref{eq: iterations} with $0<\alpha <2\ell/(d\lambda)$ results in  $f(\x_k) \convas f_*$ and $f(\x_k) \overset{L^1}{\longrightarrow} f_*$. 
\end{theorem}

\insertme {Theorem \ref{thm:convergence} guarantees $L^1$ and almost-sure convergence of the function values to a minimizer of the function whenever the function is continuously differentiable, has Lipschitz gradient, and satisfies the PL inequality (A3).  }
A broadly useful example of an objective function satisfying (A3) is linear least squares with a data matrix that is not full column rank;
Theorem \ref{thm:convergence} provides a convergence result for this \insertme{ rank-deficient linear least squares}, and similarly well-behaved non-convex functions.  Corollary \ref{corr:strong-convexity}$(ii)$ shows that the rate of convergence is linear.

\begin{corollary}[Convergence under strong-convexity and rate of convergence]\label{corr:strong-convexity}
\vspace*{-.1em}
\begin{enumerate}[label=(\roman*)]
    \item Assume (A0), (A1), (A2), (A3') and let $\x_0$ be an arbitrary initialization.
    Then recursion \eqref{eq: iterations} with $0<\alpha <2\ell/d\lambda$ results in $\x_k \convas \x_*$ where $\x_*$ is the unique minimizer of $f$.
    
    \item Assume (A0), (A1), (A2), and either (A3) or (A3').
    Then with  $\alpha=\ell/(d\lambda)$, the recursion \eqref{eq: iterations} attains the following expected rate of convergence
	\begin{equation}
	 \Expectation f(\x_k)-f_*  \leq \omega^k (f(\x_0)-f_*) ,\quad \omega = 1 - \ell \gamma/(d \lambda).
	\end{equation}

\end{enumerate}
\end{corollary}

With $\ell = d$ we recover a textbook rate of convergence, \insertme{$\omega=1-\gamma/\lambda$}, for gradient descent \cite[\S 9.3]{BoydVandenbergheBook} because, importantly, with $\ell=d$, $\P\P^\top \nabla f(\x) = \nabla f(\x)$.
\insertme{This rate is nearly optimal as can be shown with a simple example: let $d=2$ and $\x=(x,y)$ with initial conditions $\x_0=(0,1)$ and $f(\x)=\lambda/2 x^2 +\gamma/2 y^2$, then $f(\x_k) \to 0$ with linear rate $\omega=(1-\gamma/\lambda)^2$.}
Similar results to Corollary \ref{corr:strong-convexity}$(ii)$ have been derived for general stochastic gradient methods using techniques described in \cite[\S4]{NocedalBottou}. Adapting our special case to the general framework of \cite{NocedalBottou} results in the same rate of convergence as corollary \ref{corr:strong-convexity}$(ii)$; however \cite{NocedalBottou} does not address different modes of convergence, nor convergence of the iterates.
Using the more restrictive assumption of strong-convexity the result of Corollary \ref{corr:strong-convexity} is much stronger than Theorem \ref{thm:convergence}; we get almost sure convergence of the function evaluations and of the iterates to the optimal solution at a linear rate.
In inverse problems the convergence of $\x_k$, rather than that of $f(\x_k)$ is of paramount importance.
Furthermore, if either assumption (A3) or (A3') is satisfied, SSD has a linear rate of convergence. The rate of convergence is strictly better than that presented in \cite[Thm. 8]{nesterov2017random}. The rate in \cite{nesterov2017random} for $\gamma$-strongly convex objectives with $\lambda$-Lipschitz gradient is
\begin{equation}\label{eq: nesterovrate}
    \Expectation f(\x_k) - f_* \leq (\lambda/2)(1-\gamma/(8\lambda(d+4)))^k\norm{\x_0-\x_*}^2.
\end{equation}
By $\lambda$-Lipschitz gradient our Corollary \ref{corr:strong-convexity} $(ii)$ implies  
\begin{equation*}
\Expectation f(\x_k) - f_* \leq (\lambda/2)(1-\ell\gamma/(d\lambda))^k\norm{\x_0-\x_*}^2,
\end{equation*}
which is strictly better than \eqref{eq: nesterovrate}. Note that $\ell=1$ in \eqref{eq: nesterovrate}, while in our case $\ell$ can be chosen to be greater than one. 

 The proof in the convex case is different, but substantively similar to a proof of coordinate descent on convex functions found in \cite{wright2015coordinate}.
\begin{theorem}[Convergence under convexity]\label{thm: convergence-convex}
Assume (A0), (A1), (A2), (A3'').
Then recursion \eqref{eq: iterations}  with \leaveout{$0<\alpha <2\ell/(d\lambda)$} \insertme{$\alpha = \ell/(d\lambda)$} gives 

\begin{equation*}
    \Expectation f(\x_k) - f_*  \leq 2d \lambda R^2/(k \ell).
\end{equation*}

\end{theorem}
\insertme{Convergence in the convex case is in expectation, and is sub-linear. This is in line with the convergence rate of gradient descent which is also sub-linear in the smooth, convex case \cite{nesterov2013introductory}. 
In particular, taking $\ell=d$, our result gives $f(\x_k) - f_* \le 2\lambda R^2/k$ (this is now a deterministic result), where the stepsize is $\alpha = 1/\lambda$. 
It can be shown that for this common choice of a stepsize, there is a function $f$ satisfying the assumptions of Thm.~\ref{thm: convergence-convex} where $f(\x_k) - f_* \ge \frac{\lambda}{4k+2}R^2$~\cite[Thm.\ 3.2]{YoelMathProg}, which implies that when $\ell=d$, the upper bound in Thm.~\ref{thm: convergence-convex} is tight to within a factor of $8$.
}

In the general non-convex setting we can provide guarantees of convergence to a \leaveout{minimizer} \insertme{stationary point} and are able to provide guarantees on the rate at which $\norm{\nabla f(\x_k)}$ decreases. 
These are presented in the following theorem which adapts well-known results for the convergence of gradient descent on non-convex functions to our case. 
The rates of convergence are of the same order as \cite[p.24]{nesterov2017random} with slightly better constants. 
\begin{theorem}[Non-convex convergence]\label{thm: convergence-nonconvex}
Assume (A0), (A1), (A2). Then recursion \eqref{eq: iterations} with $\alpha = \ell/(d\lambda)$ and an arbitrary initialization yields
\begin{equation*}
 \min_{i\in\{0,\ldots,k\}}\Expectation \norm{\nabla f(\x_i)}^2 \leq \frac{2d\lambda (f(\x_0)-f_*)}{(k+1)\ell}.
\end{equation*}
That is, $k=\mathcal{O}(d/(\ell\epsilon))$ iterations are required to achieve $\Expectation \norm{\nabla f(\x_k)}^2 < \epsilon$.
\end{theorem}

\subsection{High-probability results}\label{sect: non-asymptotic}
While it is important to understand how an algorithm will perform on average, in practice it is good know how it is likely to perform on a single run. 
In this section we discuss convergence bounds that hold with high probability, providing a better understanding of typical convergence. We consider two types of random matrices from the class satisfying assumption (A0):
\begin{enumerate}
\item \noindent\emph{\insertme{Columns from} Haar-distributed random orthogonal matrix}:
\begin{equation}\label{eqn: Haar}
        \P = \sqrt{d/\ell}\,\mathbf{Q}\I_{d\times\ell} \in \reals^{d \times \ell},
        \end{equation}
where $\mathbf{Q}$ is as in the $QR$-decomposition of a matrix $\mathbf{Z} = \mathbf{QR} \in \reals^{d \times d}$  with $\mathbf{R}_{ii} >0,$ and each element of $\mathbf{Z}$ is drawn independently from $\mathcal{N}(0,1)$. $\I_{d \times \ell}$ truncates $\mathbf{Q}$ to its first $\ell$ columns so
$\mathbf{Q}\I_{d\times\ell}$ corresponds to $\ell$ columns of the random orthogonal matrix distributed according to the Haar measure on orthogonal matrices \cite{mezzadri2006generate}. In fact, for our results to hold, $\mathbf{R}_{ii}$ need not be strictly positive, we merely require that $\mathbf{PP}^\top \nabla f(\x) \overset{d}{=} (d/\ell)\text{Proj}_{\mathrm{col}(\mathbf{Z}\mathbf{I}_{d\times \ell})} (\nabla f(\x))$. It is convenient to work with Haar distributed matrices so we use matrices of the form \eqref{eqn: Haar}.

\insertme{There is an important correspondence between matrices described by \eqref{eqn: Haar}, Gaussian smoothing of \cite{nesterov2017random}, and the smoothing on a sphere of \cite{berahas2019theoretical}. Let  $\mathbf{Q} \in \reals^{d \times d}$ be as in \eqref{eqn: Haar}, $\v \in \reals^d$ an arbitrary fixed vector, and $\mathbf{u}=\mathbf{Z}\mathbf{e}_1$, where $\mathbf{e}_1$ is the first standard basis vector, so $\mathbf{u} \sim \mathcal{N}(\mathbf{0},\I_d)$.
Then $\mathbf{Q}^\top \v/\norm{\v}$  and $\mathbf{u}/\norm{\mathbf{u}}$ are both distributed uniformly on the $d$-dimensional sphere.
When $\ell=1$, $\P = \sqrt{d}\,\mathbf{Q}\mathbf{e}_1 \overset{d}{=} \sqrt{d}\,\mathbf{u}/\norm{\mathbf{u}}$. Therefore, in this case $(\norm{\mathbf{u}}^2/d)\P\P^\top \overset{d}{=} \mathbf{u}\mathbf{u}^\top$. That is, when $h=0$ our method is proportional to Gaussian smoothing with a constant of proportionality $\norm{\mathbf{u}}^2/d$. 
Since $\norm{\mathbf{u}}^2$ is a $\chi^2$ random variable
with mean $d$, this means in high dimensions 
the constant of proportionality is sharply concentrated around 1.}

\insertme{Furthermore, when $\ell=1$, then $\P \sim \mathcal{U}(S(0,\sqrt{d}))$, i.e., we recover spherical smoothing as discussed by \cite{berahas2019theoretical}. 
To increase $\ell$, the traditional method in the literature is to use \eqref{eq: berahas gaussian smoothing}, which is different than \eqref{eqn: Haar} for $\ell>1$.
Thus, we provide a novel generalization of Gaussian smoothing and smoothing on the sphere as a mapping of the gradient onto a lower dimensional subspace.
A consequence of our approach is that 
matrices of the form \eqref{eqn: Haar} with $\ell=d$ 
satisfy 
$\P\P^\top=\P^\top\P=\I_d$, and the exact gradient is recovered.}

\insertme{
To summarize, in high-dimensions, our method with matrices defined by \eqref{eqn: Haar} is very similar to both Gaussian smoothing and smoothing on a sphere for $\ell=1$, but the differences with existing methods grow as $\ell$ increases, and are markedly different for $\ell=d$.
Table \ref{table: special cases} provides a summary of the well-known special cases to our algorithm, and describes how they relate to our framework. 
We re-emphasize that for both Gaussian smoothing and smoothing on a sphere it is typical to use $h>0$, however we only analyze the case $h=0$. 
Indeed it is true the $h=0$ is no longer smoothing of the gradient \emph{per se}, but is a projection of the gradient onto a subspace of dimension $\ell$.}

\begin{table}[H]
\centering
\begin{tabular}{llll}
Description & $\ell$ & $\P_k$ & $\alpha$ \\ \hline
Gaussian smoothing & 1 & Satisfying \eqref{eqn: Haar} & $\norm{\mathbf{Z}_k\mathbf{e}_1}^2/(d^2 \lambda)$ \\
Smoothing on a sphere & 1 & Satisfying \eqref{eqn: Haar} &  $0<\alpha<2/(d\lambda)$ \\  \, \, \, of radius $\sqrt{d}$ \\
Gradient descent & d & Any satisfying (A0) & $0<\alpha<2/\lambda$ \\
Block-coordinate descent & $1\leq \ell<d$ & Satisfying \eqref{eqn: CD} & $0<\alpha<2\ell/(d\lambda)$
\end{tabular}

  \captionsetup{justification=centering}
 \caption{Summary of special cases of our framework. Using the $\ell$, $\P_k$, and $\alpha$ specified in the table it is possible to recover exactly the methods described. $\mathbf{Z}_k$ is the $k^{\text{th}}$ Gaussian matrix used to generate $\P_k$ and $\mathbf{e}_1$ is the first standard basis vector. \label{table: special cases} }
\end{table}

Note that for problems of interest, function evaluations are so costly that we can ignore the computational overhead of a QR decomposition,  which is $\mathcal{O}(d\ell^2 - 2\ell^3/3)$. Since $\ell\ll d$, the cost is negligible compared to, for instance, $d$ PDE-solves.

\item \noindent\emph{Randomized block-coordinate descent random matrix}:
    \begin{equation} \label{eqn: CD}
        \P = \sqrt{d/\ell}\, \mathbf{D},
    \end{equation}    
    where $\mathbf{D} \in \reals^{d \times \ell}$ is comprised of $\ell$ columns of the identity matrix $\I_d$ selected uniformly at random. 
It is straightforward to verify that \eqref{eqn: Haar} and \eqref{eqn: CD} satisfy assumption (A0), the former by properties of the QR decomposition.
\medskip
Denoting the columns of $\P$ as $\P^1, \ldots, \P^{\ell}$, the following equality holds

\begin{equation}\label{eq: finite difference P}
   \nabla f(\x) \approx \P\P^\top \nabla f(\x) = \begin{pmatrix} \P \nabla_{\P^1} f(\x) \\
   \vdots \\
   \P \nabla_{\P^{\ell}} f(\x)
   \end{pmatrix},
   \end{equation}
\insertme{where $\nabla_{\P^i}f(\x)$ is a directional derivative of $f$ at $\x$ in the $\P^i$-th direction}. Thus there is a convenient interpretation that the gradient is approximated by a mapping onto an $\ell$-dimensional random subspace embedded in $\reals^d$. In fact, since $\Expectation \P\P^\top = \I$, $\P\P^\top \nabla f(\x)$ is \leaveout{an unbiased estimator of} \insertme{centered at} $\nabla f(\x)$ \leaveout{whose MSE is} \insertme{with MSE} $(1-\ell/d) \norm{\nabla f(\x)}^2$.
\leaveout{In particular when $h=0$, Gaussian smoothing given by \eqref{eq: gaussiansmoothing fd} is identical to mapping the gradient onto a 1-dimensional random orthogonal subspace of $\reals^d$ via a Haar-distributed random matrix.
Choosing $\ell>1$ and using matrices of the form \eqref{eqn: Haar} generalizes the notion of Gaussian smoothing such that when $\ell=d,\,\P\P^\top=\P^\top\P=\I$ and the exact gradient is recovered.}
\end{enumerate}

In advance of the main results of this section we investigate how well multiplication by the matrices specified by \eqref{eqn: Haar} and \eqref{eqn: CD} preserves the norm of an arbitrary vector. 
Norm invariance has important consequences with respect to the rate of convergence. Of particular interest for our purpose is the lower bound which governs the rate of convergence (see Theorem \ref{thm: strong-convexity} for details). We define a successful embedding in order to quantify the norm invariance.

\begin{definition}[Successful isometric embedding] An embedding $\P$ is deemed to be successful if for some $\epsilon \in (0,1)$ and some $\v \in \reals^d$, $\norm{\P^\top \v}^2 \geq (1-\epsilon) \norm{\v}^2$.
\end{definition}

The following Lemma provides the probability of successful embedding when the matrix $\P$ is Haar-distributed.
\begin{lemma}[Approximately isometric embedding using Haar-distributed matrices]\label{Lemma: DavidJL}
Fix $\epsilon \in (0,1)$, \leaveout{an integer $\ell \in [1,d]$} \insertme{a positive integer $\ell \leq d$}, and consider a matrix $\mathbf{P}$ drawn according to \eqref{eqn: Haar}.
Then
for any fixed vector $\v \in \reals^d$, the probability of a successful embedding,  $\delta$,  
is given by
\begin{equation*}
  \delta =  1-I_{(1-\epsilon)\ell/d}(\ell/2,(d-\ell)/2)= \Prob(X \geq (1-\epsilon)\ell/d),
  \end{equation*}
 where $I_p(\alpha,\beta)$ is the regularized incomplete Beta function, and $X \sim \mathcal{B}eta(\ell/2, (d-\ell)/2)$. 
\end{lemma}

\begin{figure}[ht]
      \centering
      \includegraphics[width=.8\linewidth]{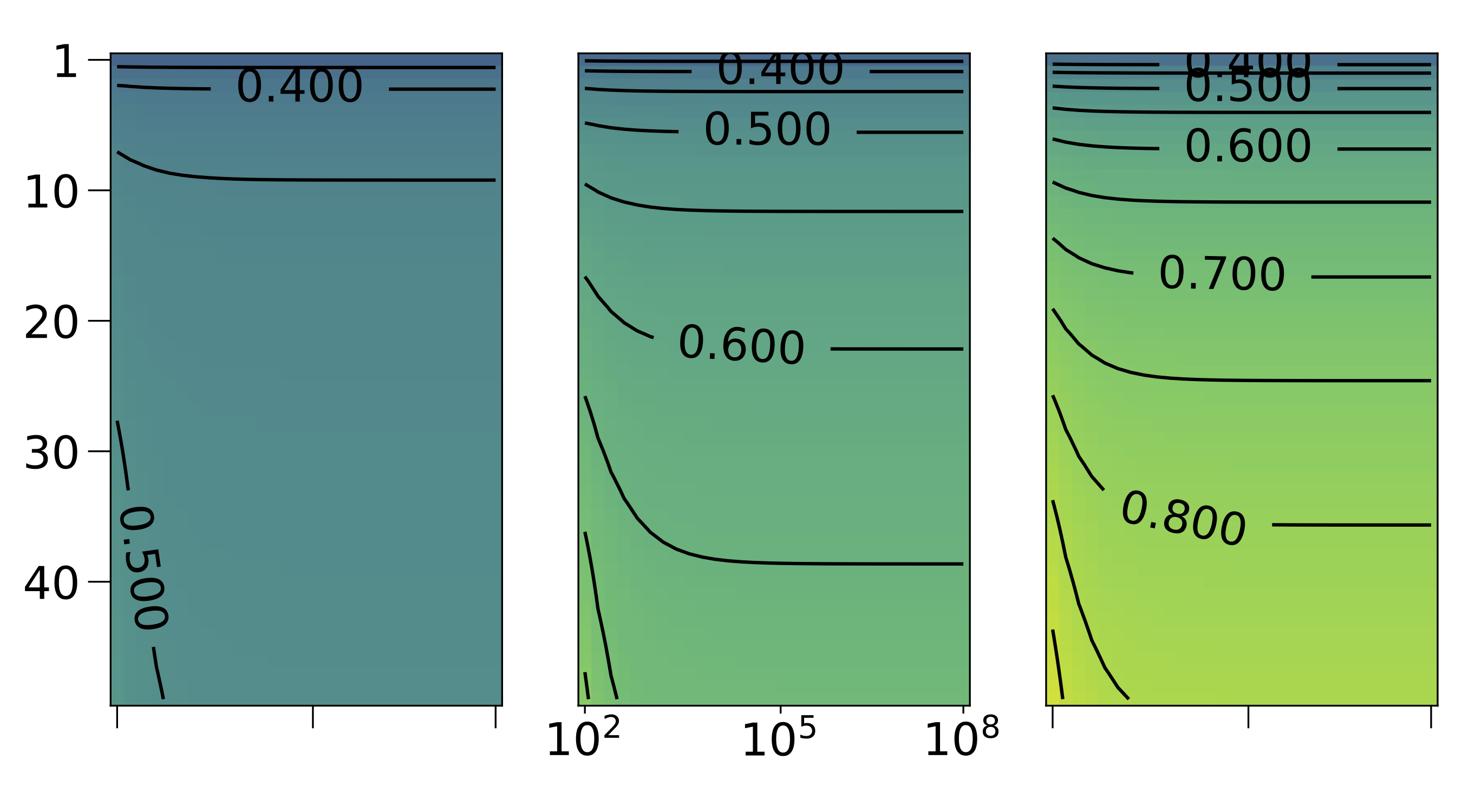}
        \put(-145,0){\small $d$}
        \put(-285,88){\small $\ell$}
      \vspace*{-.5em}
      \caption{Contour plots for probability of successful embedding for various values of $\ell$, $d$, and $\epsilon$. Each of the figures share the same horizontal and vertical range. 
      \textbf{Left}: $\epsilon=0.01$. 
      \textbf{Center}: $\epsilon=0.1$.
      \textbf{Right}: $\epsilon=0.2$. }
      \label{fig: embeddings}
      \vspace{-.5em}
  \end{figure}

For fixed $d$ one can simply use Lemma \ref{Lemma: DavidJL} to determine values of $\ell$ and $\epsilon$ required to achieve the desired probability of successful embedding. \insertme{For a fixed $d$, an increase in $\ell$ or $\epsilon$ corresponds to an increase in $\delta$.
For $\ell=d$, we can take $\epsilon=0$ and $I_1(\alpha,\beta)=0$ so $\delta=1$, meaning we always have a perfect isometric embedding.
} Figure \ref{fig: embeddings} provides examples of the probability of success for various values of $\ell$, $d$, and $\epsilon$. It is plain to see the similarity between the left hand-side of Lemma \ref{Lemma: DavidJL} and the lower tail of the Johnson-Lindenstrauss (JL) lemma when it is applied to a single point. Indeed, a connection of the JL lemma with the Beta distribution is discussed in \cite{frankl1990some}. Our bound differs in two ways: first, in \cite{frankl1990some} they provide asymptotic results as $d \to \infty$ whereas our results are valid for all $d$ with the $d$-dependence explicit; second, \cite{frankl1990some} provide a closed-form approximate bound while we provide an exact functional form. For finite dimensions, our result is stronger.

A well-known property of the matrices \eqref{eqn: Haar} is that $\P^\top$ is spherically symmetric.
That is $\P^\top \mathbf{U}$ has the same distribution as $\P^\top$ for any orthogonal matrix $\mathbf{U}$.
Consequently, the quality of the embedding does not depend on the vector $\v \in \reals^d$.
Naturally, the coordinate descent matrices given by \eqref{eqn: CD} do not share this orthogonal invariance; indeed, speaking of the ability of such matrices to preserve pairwise distances, Achlioptas \cite{achlioptas2003database} says ``A naive, perhaps, attempt at constructing JL-embeddings would be to pick $\ell$ of the original coordinates in $d$-dimensional space as the new coordinates.
Naturally, as two points can be very far apart while only differing along a single dimension, this approach is doomed''. 
Remark \ref{remark: CD} provides intuition for the reason randomized block-coordinate descent cannot be close to norm preserving for arbitrary directions. 

\begin{remark}[Coordinate sampling is rarely an isometry]\label{remark: CD}
 Let $\mathbf{v} \in \reals^d$ be a standard basis vector and $\P \in \reals^{d \times \ell}$ 
 be a coordinate descent sampling matrix satisfying \eqref{eqn: CD}.
Then, $\norm{\P^\top \v} \in \{0,1\}$, and
 \begin{equation*}
      \Prob(\|\sqrt{\ell/d}\,\P^\top\mathbf{v}\|^2=1) = \ell/d \quad \text{and} \quad \Prob(\|\sqrt{\ell/d}\,\P^\top\mathbf{v}\|^2=0) = 1-\ell/d.
 \end{equation*}
 Thus, 
 \begin{equation*}
     \Expectation\|\P^\top \v\|^2=1 \quad \text{and} \quad \mathbb{V}\mathrm{ar}\|\P^\top\v\|^2 = d/\ell-1.
 \end{equation*}
\end{remark}
Since exactly $\ell$ entries of $\mathbf{\P\P^\top}$ are $1$, the probability that any non-zero entry corresponds to a non-zero entry of $\mathbf{v}$ is $\ell/d$.
\bigskip

Remark \ref{remark: CD} shows that in the worst case (that is, if the vector $\mathbf{v}$ is axis-aligned with concentration along a single coordinate), there is no approximate norm-preservation: it is either exact with probability $\ell/d$ or not-at-all with probability $1-\ell/d$. This compares very unfavorably to the results of Lemma \ref{Lemma: DavidJL}, \leaveout{as shown in} cf. Figure \ref{fig: embeddings}.
\medskip

To summarize, the structure of the objective function plays a role in the quality of a coordinate descent mapping, and in the worst-case the mapping using \eqref{eqn: CD} is useless with probability $1-\ell/d$. 
In contrast, using Haar matrices guarantees that irrespective of the structure of the function a successful embedding is obtained with probability according to Lemma \ref{Lemma: DavidJL}\leaveout{, which is not sensitive to dimension}. \insertme{Though this probability depends on dimension, it is not very sensitive to an increase in $d$, as illustrated in Figure \ref{fig: embeddings}}.
\smallskip

Due to  the strong dependence on the dimension for randomized coordinate descent, the analysis in the remainder of this section is not appropriate for matrices of type \eqref{eqn: CD}. Thus, we consider only Haar distributed random matrices. It should be noted that there are special classes of functions for which the complexity is independent of $d$, as discussed in \cite{BCGDIterationCompliexity}, however in general the dependence on the dimension can not be removed using coordinate descent methods. 
We consider first a result that is a simple but useful corollary to Theorem 3.1 in \cite{marchal2017sub}, later proved in \cite{arbel2020strict}

\begin{remark}\label{Lemma: Optimal binomial}
Let $B \sim Bin(k,\delta)$. Then for all $t > 0$ and $\delta \in (0,1)$
\begin{equation*}
    \mathbb{P}(B > k\delta+t) \leq \exp\left(-t^2/(2\sigma_k^2)\right)
\quad \text{and}\quad    \mathbb{P}(B < k\delta-t) \leq \exp\left(-t^2/(2\sigma_k^2)\right)    
\end{equation*}
with
\begin{equation}\label{eq: sigma}
    \sigma_k^2 = \begin{cases} \frac{k (1-2\delta)}{2 \log ((1-\delta)/\delta)} & \delta \in (0,1)\setminus \{1/2\} \\
    1/4, & \delta = 1/2.
    \end{cases}
\end{equation}
\end{remark}

Remark \ref{Lemma: Optimal binomial} provides an optimal proxy-variance for sub-Gaussianity of Binomial random variables. 
For $\delta=1/2$, $\sigma_k^2$ is defined as $k/4$ so that $\sigma_k$ is continuous in $\delta$;
\insertme{
also note that for any $k$, $\lim_{\delta\nearrow 1}\sigma_k^2=0$ which agrees with the fact $\mathbb{P}(B\neq k)=0$ in the case $\delta=1$ (which occurs when $\ell=d$). 
}
We use the result of Remark \ref{Lemma: Optimal binomial} to provide sharp bounds for the performance of our algorithm.
First we state a result showing that the success of each embedding is independent so that the number of successful embeddings can be treated as a binomial random variable, which in turn allows for an application of Remark \ref{Lemma: Optimal binomial}.

\begin{remark}\label{remark: independence}
Let $A_k(\v_k) = \left\{\norm{\P_k^\top \v_k}^2 \leq (1-\epsilon) \norm{\v_k}^2 \right\}$ and $\v_k$ be independent of $\P_k$ for all $k$ with $\P_k$ drawn according to \eqref{eqn: Haar}. Then $(A_k(\v_k))$ is an independent sequence of events.
\end{remark}

The remark is proved by iteratively conditioning on the available information and recognizing that spherical symmetry implies $A_k$ is identically distributed for any $\v_k$ that is fixed or independent of $\P_k$\insertme{, and can be found in the Appendix}. Using Lemma \ref{Lemma: DavidJL} and Remark \ref{Lemma: Optimal binomial} in conjunction with Remark \ref{remark: independence} results in the following probabilistic rate of convergence,

\begin{theorem}[Probabilistic rate of convergence. Strongly-convex case]\label{thm: strong-convexity}
     Assume (A1), (A2), (A3') and let $\x_0$ be an arbitrary initialization. Apply recursion \eqref{eq: iterations} with step-size
    $\alpha=\ell/(d\lambda)$ and $\P_k$ drawn according to \eqref{eqn: Haar}, with $\ell$ sufficiently large to achieve the desired $\epsilon$ and $\delta$  according to Lemma $\ref{Lemma: DavidJL}$. Then for any $t \in (0,\delta]$
\begin{equation*}
\Prob\left(f_e(\x_{k}) \geq \rho^k f_e(\x_0) \right) \leq \exp(-(kt)^2/2\sigma_k^2),
\end{equation*}
where $\sigma_k^2$ is defined by \eqref{eq: sigma} and,
\begin{equation*}
 \rho = \left(1-(1-\epsilon)\frac{\ell \gamma}{d \lambda}\right)^{\delta  - t}.
\end{equation*}
\end{theorem}

Theorem \ref{thm: strong-convexity} provides an exponential decay (in $k$) for the probability that any single run of the algorithm converges more slowly than the average performance guaranteed by Corollary \ref{corr:strong-convexity}$(ii)$. Similar results can be derived for the convex case combining the methodology of Theorem \ref{thm: strong-convexity} with Theorem \ref{thm: convergence-convex}.
\insertme{We note the interplay between parameters $\delta, t$ and $\epsilon$, all of which affect $\rho$: we can trade off $\epsilon\to 0$ by decreasing $\delta$; both $\epsilon$ and $\delta$ affect $\rho$, but due to their complicated relationship, it is not easy to optimize $\rho$ with respect to these parameters. We can also take $t\to \delta$ to get a conservative bound (high probability but worse rate $\rho$), or $t\to 0$ to get an aggressive bound (lower probability but better rate $\rho$). Since the probability concentrates quickly with the number of iterations $k$, for large iterations, one can take $t\propto 1/\sqrt{k}$ and have both good control over the failure probability, $\text{exp}(O(-k))$, while still having a good convergence rate.}

\insertme{Again, Theorem~\ref{thm: strong-convexity} agrees with the standard deterministic result (discussed after Corollary~\ref{corr:strong-convexity}) when $\ell=d$, since then $\epsilon=0$, $\delta=1$ and  $\sigma_k^2=0$ for any $t>0$, so the rate $\rho$ is arbitrarily close to $1-\frac{\gamma}{\lambda}$, 
which is the rate from Corollary~\ref{corr:strong-convexity}$(ii$).
}

\section{Experimental results}\label{sect: EmpiricalResults} In this section we provide results for a synthetic problem, a problem from the machine learning literature, and a PDE-constrained shape-optimization problem.
In the synthetic and machine learning problems we compare to randomized block-coordinate descent.
For the shape-optimization problem we compare to Gaussian smoothing and finite-difference gradient descent. \insertme {In each of the examples we make use of a deterministic backtracking line search with Armijo conditions. Analysis of algorithms with inexact gradient estimates using a stochastic line search is a topic that has received considerable attention recently, but which we do not address here. In particular, \cite{berahas2019global,PaquetteScheinberg} describe a stochastic variant of an Armijo backtracking line search that can be adapted to our method to provide sharper convergence analysis. Their work does not directly apply to all of our settings without modification, but in the strongly-convex case the application is clear.}
\subsection{Synthetic data}\label{subsect:synthetic-data}
\begin{figure}[ht]
      \centering
      \includegraphics[width=.32\linewidth]{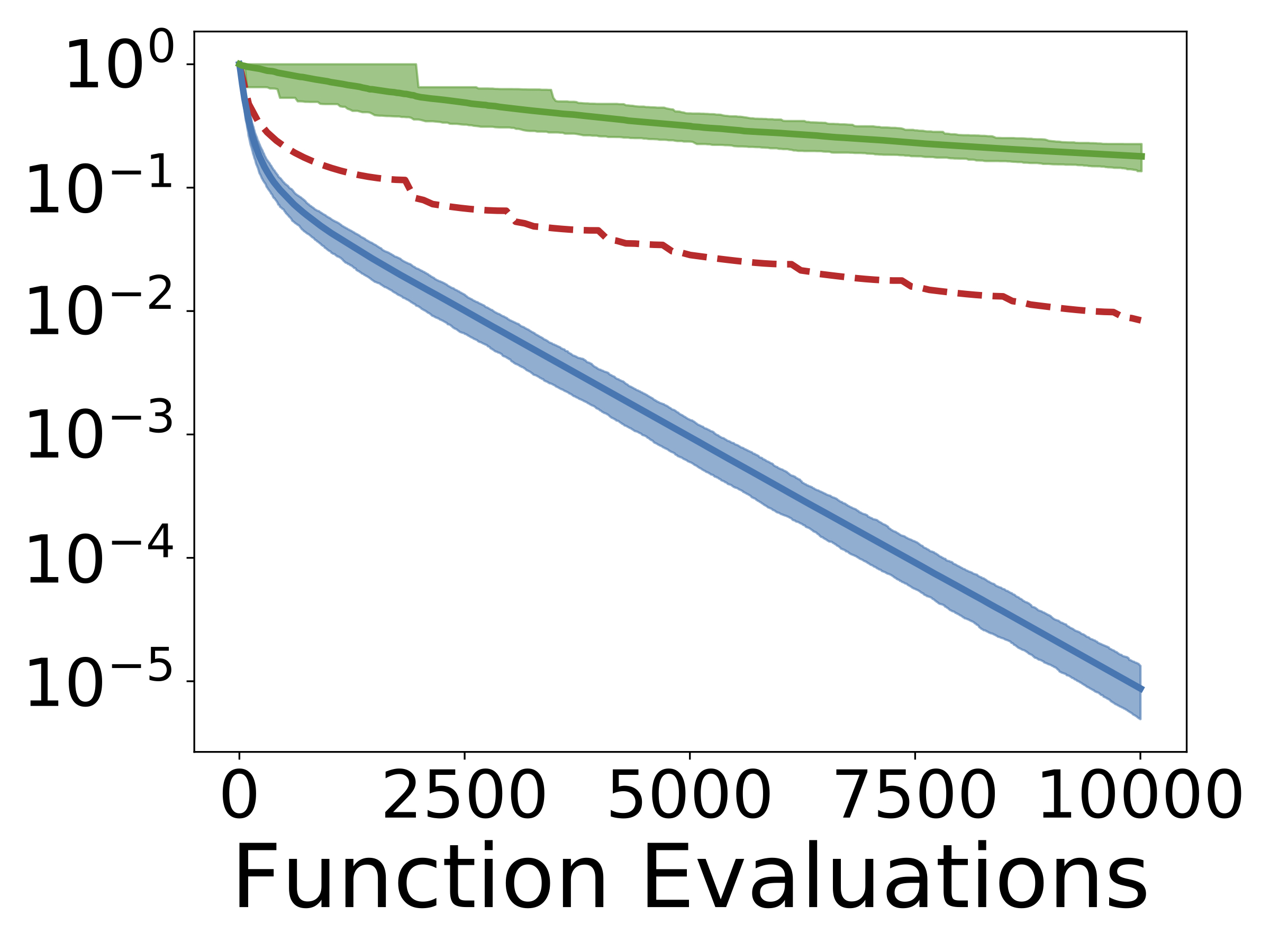}
      \includegraphics[width=.32\linewidth]{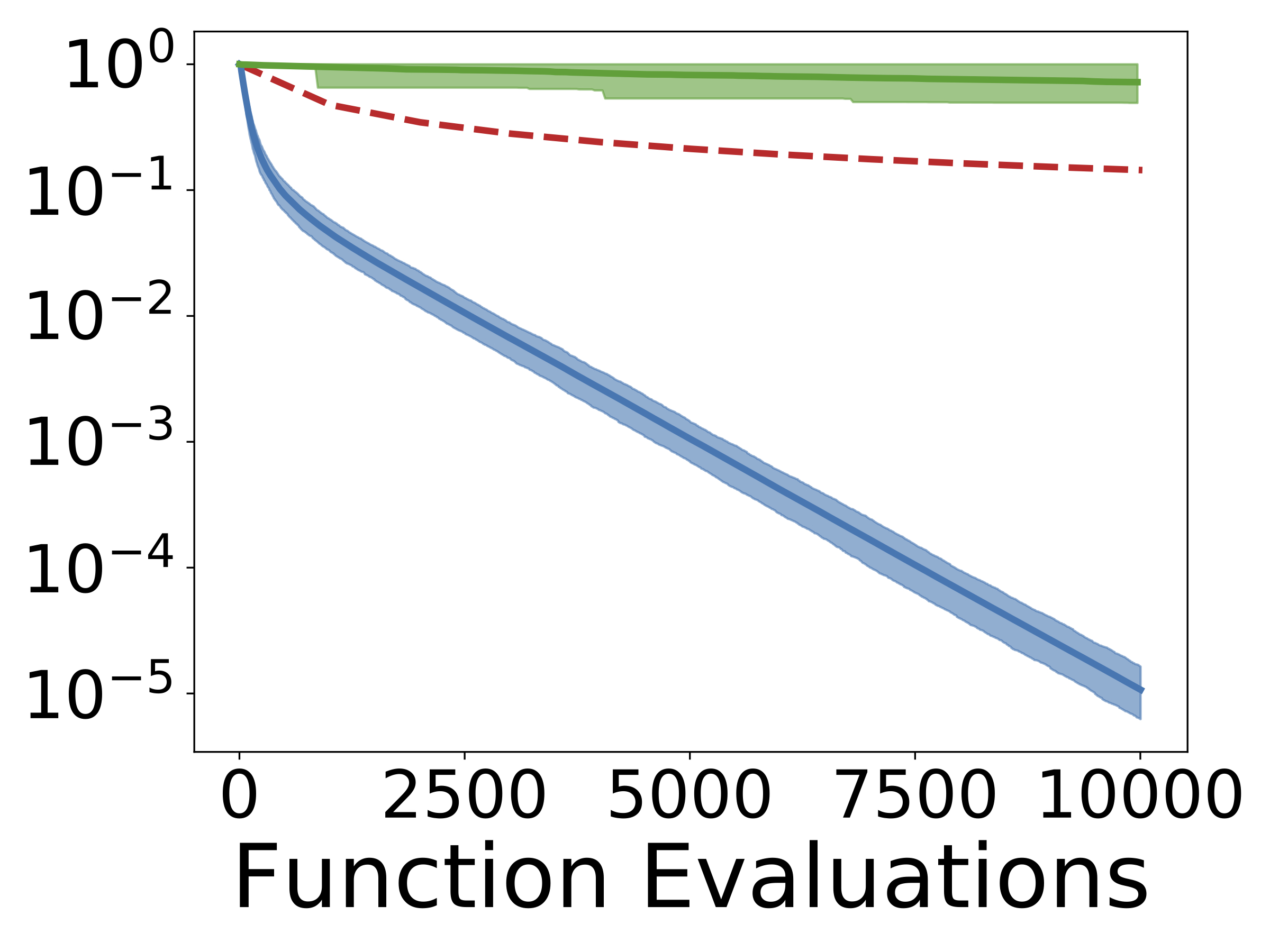}
      \includegraphics[width=.32\linewidth]{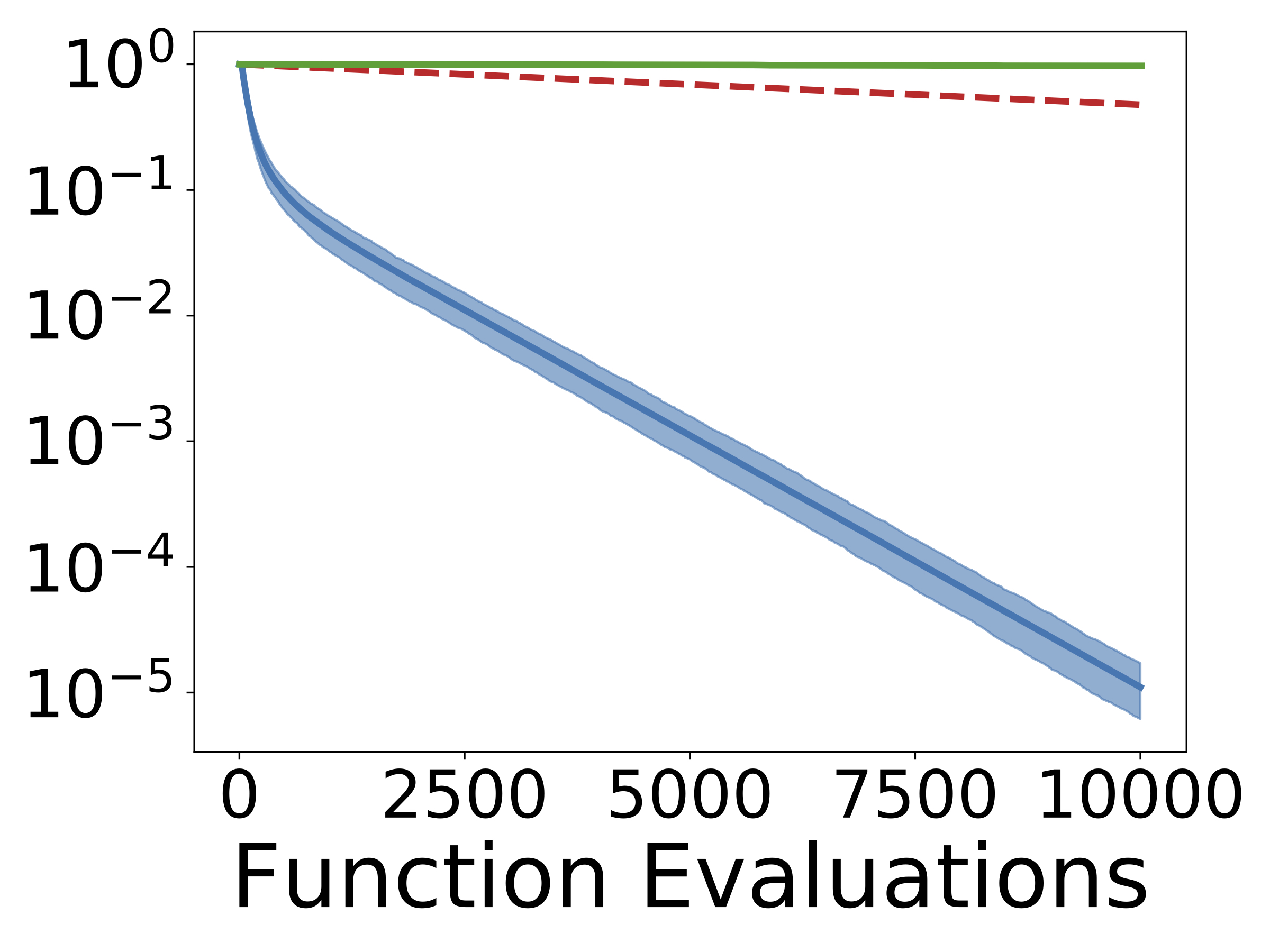}
            \includegraphics[width=.32\linewidth]{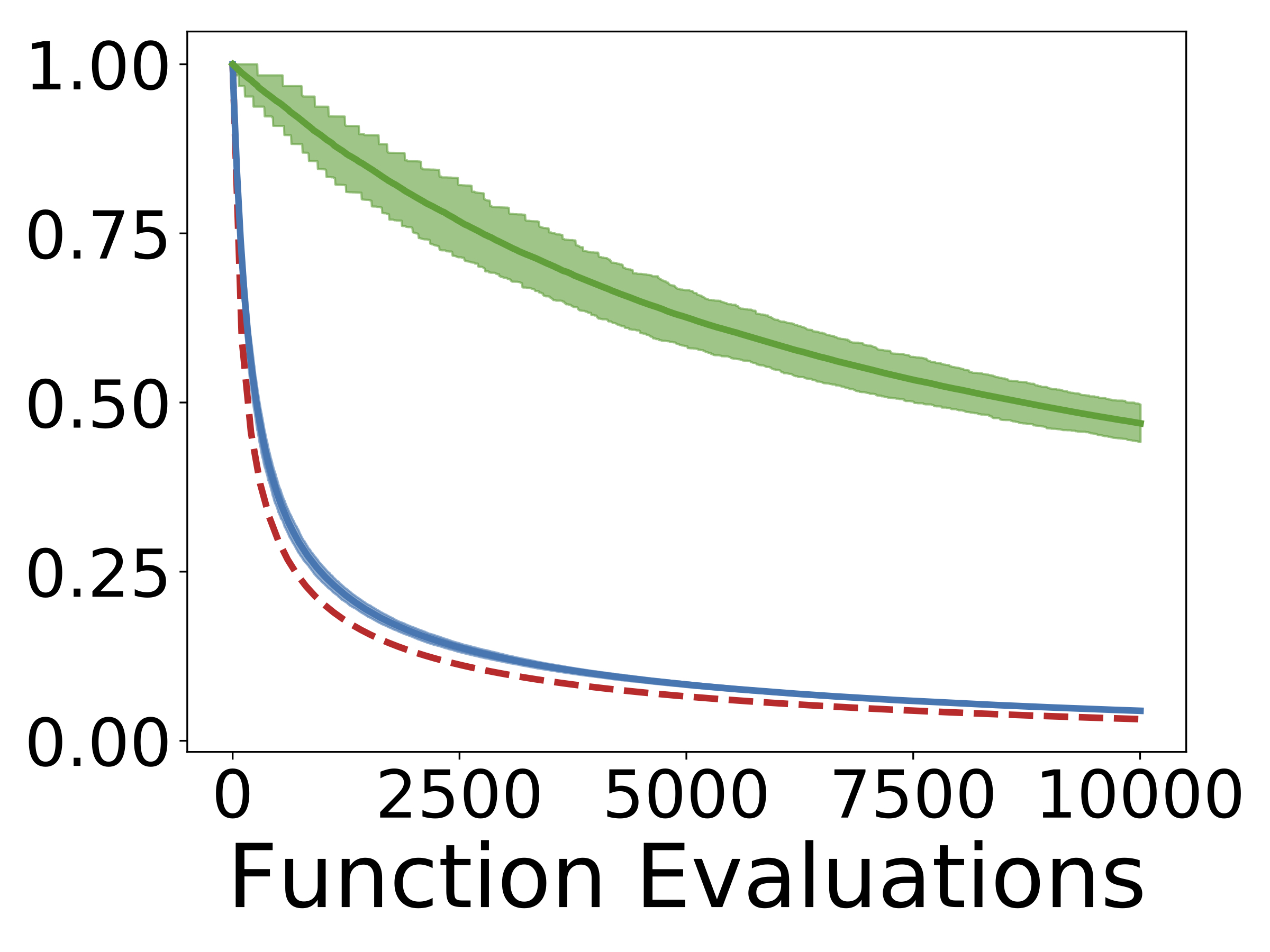}
      \includegraphics[width=.32\linewidth]{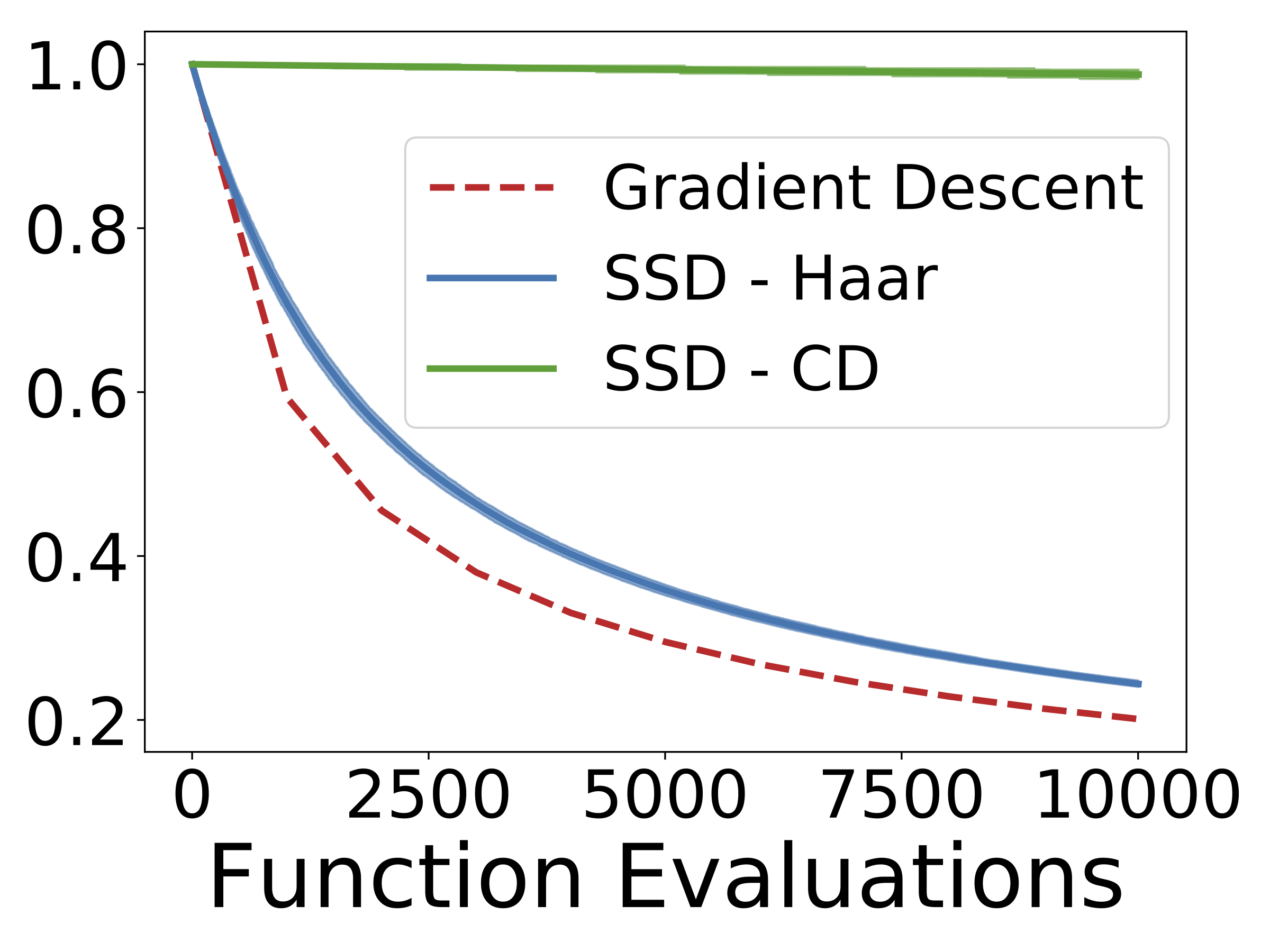}
      \includegraphics[width=.32\linewidth]{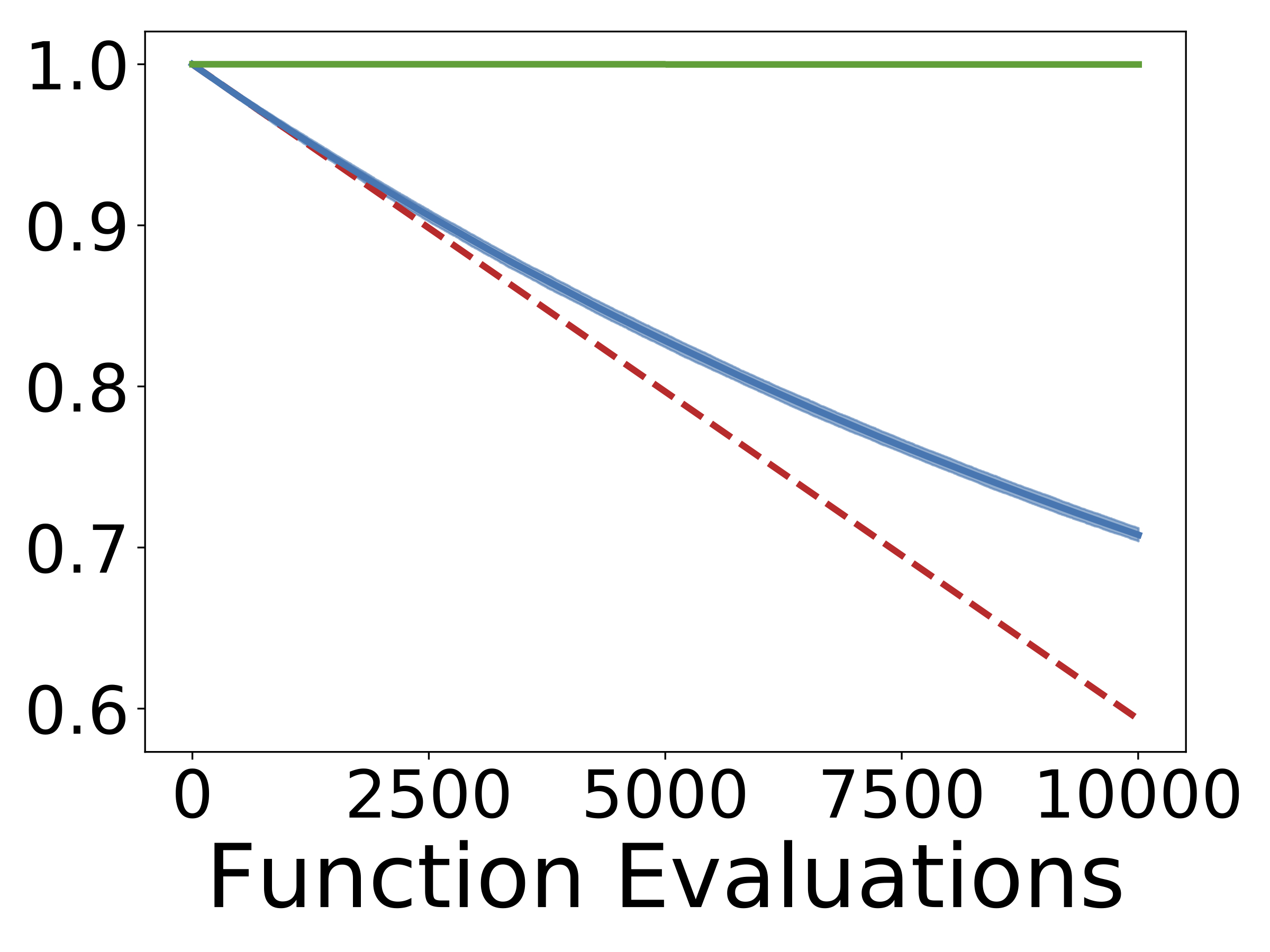}
      \vspace*{-1em}
      \caption{
      Minimizing a function from the family \eqref{eq: NesterovWorst} with $r=20,~\lambda =8$. CD represents randomized block-coordinate descent. In several of the subfigures gradient descent overlaps randomized block-coordinate descent. The shaded regions in the SSD cases represent the interval between best  $10^{\text{th}}$ and $90^{\text{th}}$ percentile performance after $1000$ runs. The vertical-axis is the relative error: $(f(\x_k)-f_*)/f_*$.
      \textbf{Left}: $d=100$.
      \textbf{Center}: $d=1000$.
      \textbf{Right}: $d=10000$.
      \textbf{Top}: Step-size chosen by a backtracking linesearch with Armijo conditions.
      \textbf{Bottom}: Fixed step-size. }
      \label{fig: nesterov}
  \end{figure}

We begin with a simulated example using what Nesterov dubs "the worst function in the world" \cite{nesterov2013introductory}.  Fix a Lipschitz constant $\lambda>0$ and let
\begin{equation}\label{eq: NesterovWorst}
    f_{\lambda,r}(\x) = \lambda (( x_1^2 + \medmath{\sum_{i=1}^{r-1}} \left( x_i - x_{i+1}\right)^2 + x_r^2 )/2 - x_1 )/4,
\end{equation}
where $x_i$ represents the $i^{\text{th}}$ coordinate of $\x$ and $r < d$ is a constant integer that defines the intrinsic dimension of the problem. This function is convex and continuously differentiable with global minimum $f_*= -\lambda r/8(r+1)$, so Theorem \ref{thm: convergence-convex} applies. 
This example illustrates the consequences of the dimension dependence in Remark \ref{remark: CD}, as well as the dimension independence of Lemma~\ref{Lemma: DavidJL} in the context of optimization using recursion \eqref{eq: iterations}. Figure \ref{fig: nesterov} highlights the performance of three algorithms: finite-difference gradient descent, SSD using \eqref{eqn: CD} (hereafter, SSD-CD), and SSD using \eqref{eqn: Haar} (hereafter, SSD-Haar); all algorithms start at $\x=0$. We show each with the fixed step-size $\alpha=\ell/(d\lambda)$ suggested by the theorem, as well as an adaptive step-size using a backtracking linesearch with the Armijo conditions. We keep $\ell=3$ and $r = 20$ fixed and provide results for $d=100$, $d=1000$, $d=10000$. For the SSD cases we run each 500 times and display the performance of the $10^{\text{th}}$ and $90^{\text{th}}$ percentile (shaded region) as well as the mean performance.

Clearly both gradient descent and randomized block-coordinate descent depend strongly on the ambient dimension of the problem, even when a linesearch is used. 
Functions from this family are a worst case for both of these algorithms as only the first $r$ dimensions have a non-zero gradient.
Thus, in the case $d=10000$, gradient descent must perform $10000$ function evaluations at every iteration when only $r=20$ dimensions are important.
Similarly, randomized coordinate descent has only a 20/10000 chance of descending at all, so as predicted in the discussion of Remark \ref{remark: CD}, we see many iterations of coordinate descent with no improvement. The linesearch makes coordinate descent slower relative to gradient descent for this example because every iteration for which a pertinent coordinate is not selected requires several function evaluations to perform the linesearch. 
Regarding SSD-Haar, using a linesearch dramatically impacts performance by allowing for invariance to ambient dimension as suggested by Lemma \ref{Lemma: DavidJL}. Without linesearch, as expected by Theorem \ref{thm: convergence-convex}, the performance can be no better than that of gradient descent.
As previously noted, the function has low intrinsic dimension; the performance on this problem suggests that the bound in Lemma \ref{Lemma: DavidJL} (and in turn, of Theorems \ref{thm:convergence} and \ref{thm: convergence-convex}) can be sharpened by accounting for this structure and we consider this a promising avenue for future research.

\subsection{Parameter estimation for sparse Gaussian processes}\label{subsect: experiments-GP}
We test the efficacy of SSD-Haar against SSD-CD in the context of hyperparameter estimation for sparse Gaussian processes used in regression.
The goal is inference on a function $T:\reals^d\to\reals$  based on noisy observations at $m$ points
$\z_1,\ldots,\z_m$.  We use a zero-mean Gaussian process with covariance function $
    \mathbb{C}\mathrm{ov}(T(\z_i),T(\z_j)) = K(\z_i, \z_j; \boldsymbol{\theta})$ and model the $m$ observations as $y_i = T(\z_i) + \epsilon_i$, where $K(\cdot\,,\cdot\,;\boldsymbol{\theta})$ is a symmetric positive-definite kernel with parameters $\boldsymbol{\theta}$.
The process $T$ is assumed to be independent of the noise vector $(\epsilon_1,\ldots,\epsilon_m)^\top \sim N(\zero,\sigma^2\Ib)$ with unknown variance $\sigma^2$.
We denote the covariance of the vector $\mathbf{T}=(T(\z_1),\ldots,T(\z_m))^\top $ as $\boldsymbol{\Sigma}_{\mathbf{T}} =  \mathbb{V}\mathrm{ar}(\mathbf{T})$, where $(\boldsymbol{\Sigma}_{\mathbf{T}})_{ij} = K(\z_i, \z_j;\boldsymbol{\theta})$.
Maximum likelihood estimates of the parameters $\boldsymbol{\Theta} = [\boldsymbol{\theta}, \sigma^2]$ are obtained by maximizing the 
log-marginal likelihood of observations 
$\y=(y_1,\ldots,y_m)^\top $ with density $p_{\y}$ \cite{Rasmussen2009GaussianPF}: $\ell(\boldsymbol{\Theta};\y) = \log p_{\y}(\y;\boldsymbol{\Theta})$. %
When the number of observations is large the cost of this maximization is $\mathcal{O}(m^3)$ due to the inversion and determinant calculations in $\ell(\boldsymbol{\Theta};\y)$. We use the method described in \cite{titsias2009variational} to approximate the likelihood. The basic idea is as follows: choose a $p < m$ and define a set of inducing points $\widetilde{\z}_1,\ldots,\widetilde{\z}_p \in \reals^d$ different from the original $\z_1,\ldots,\z_m$, and let
$\widetilde{\mathbf{T}} = (T(\widetilde{\z}_1),\ldots, T(\widetilde{\z}_p))^\top $. We obtain a lower bound for the loglikelihood \cite{titsias2009variational}:
\begin{equation}\label{eq:llbound}
   \ell(\boldsymbol{\Theta};\y) \geq f(\widetilde{\z}_1,\ldots\widetilde{\z}_m,\boldsymbol{\Theta})=\widetilde{\ell}(\boldsymbol{\Theta};\y) - \mathrm{tr}
   (\mathbb{V}\mathrm{ar}(\mathbf{T}\mid \widetilde{\mathbf{T}}))/2\sigma^2.
\end{equation}
Here
$\widetilde{\ell}$ is the loglikelihood of the multivariate Gaussian $N(\zero, \widehat{\boldsymbol{\Sigma}}_{\mathbf{T}})$, where $
\widehat{\Sigmab}_{\mathbf{T}} = \Sigmab_{\mathbf{T}} - \mathbb{V}\mathrm{ar}(\mathbf{T}\mid \widetilde{\mathbf{T}})=
\mathbb{C}\mathrm{ov}(\mathbf{T},\widetilde{\mathbf{T}}) \,\Sigmab_{\widetilde{\mathbf{T}}}^{-1}\,\mathbb{C}\mathrm{ov}(\widetilde{\mathbf{T}},\mathbf{T})
$ 
is the  the Nystr\"om approximation of $\Sigmab_{\fb}$ introduced in \cite{williams2001using}. 
Gradient-based methods are used to simultaneously find an optimal placement of the $p$ inducing points and the best hyperparameter settings by maximizing the lower
bound in \eqref{eq:llbound}, which we re-state as a function of $\x = [\widetilde{\z}_1, \ldots, \widetilde{\z}_p, \boldsymbol{\Theta}]$ to be consistent with notation in previous sections:
\begin{equation}\label{eqn: VariationalUpperBound}
f(\x) = \widetilde{\ell}(\boldsymbol{\Theta};\y) -\,\mathrm{tr}
   (\mathbb{V}\mathrm{ar}(\mathbf{T}\mid \widetilde{\mathbf{T}}))/2\sigma^2.
\end{equation}
 Practically speaking, the optimization problem is $(pd+\abs{\boldsymbol{\theta}}+1)$-dimensional: $pd$ for $p$ inducing points in $\reals^d$, $\abs{\boldsymbol{\theta}}$ for the kernel hyperparmeters, and 1 for the unknown noise variance. 
 By moving to this high-dimensional optimization problem the time complexity is reduced to $\mathcal{O}(mp^2)$ and the storage costs to $\mathcal{O}(mp)$. 

 \medskip
For example, we model a noisy version of the function described by \eqref{eq: NesterovWorst} with $\lambda = 1$ and $r=d$ using a Gaussian process in the framework of \cite{titsias2009variational} with a squared-exponential kernel that has two unknown parameters. Between the inducing points, the parameters of the kernel, and the unknown noise, there are $153,~503,~2003$ parameters to be estimated for cases $(d=3,~p=50),~(d=10,~p=50)~,(d=20,~p=100)$ respectively. We report the objective function, which is \eqref{eqn: VariationalUpperBound} up to an irrelevant constant. \leaveout{We perform the optimization for 500 function evaluations, so would not have the opportunity to take a single step in the second and third experiment} \insertme{We terminate the algorithm after 500 function evaluations. Thus, since the second and third experiments have 503, and 2003 parameters respectively, gradient descent would not have the opportunity for even one iteration.} As such, for all three experiments we only compare SSD-Haar to SSD-CD.
 
    \begin{figure}[ht]
      \centering
      \includegraphics[width=.32\linewidth]{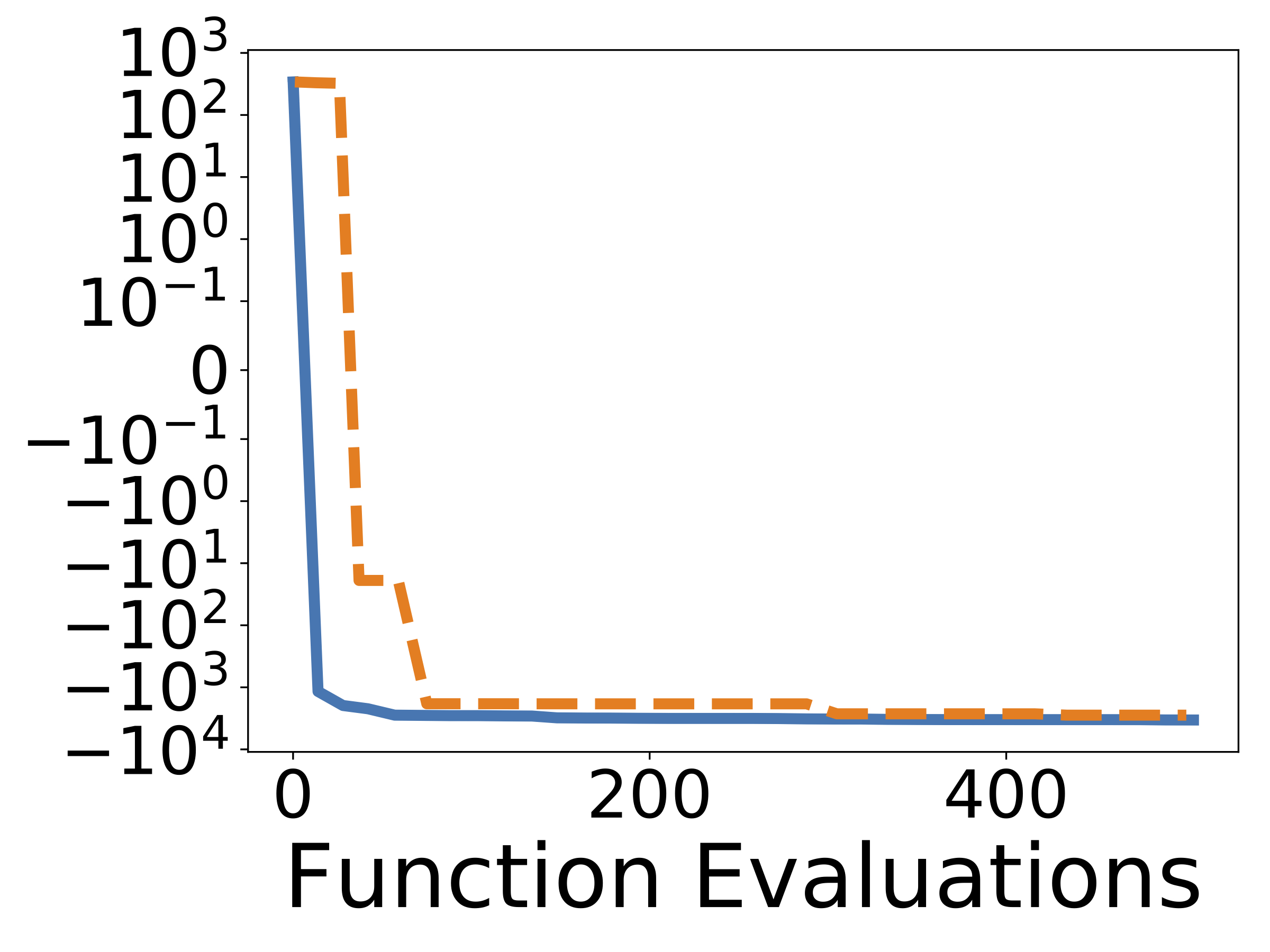}
      \includegraphics[width=.32\linewidth]{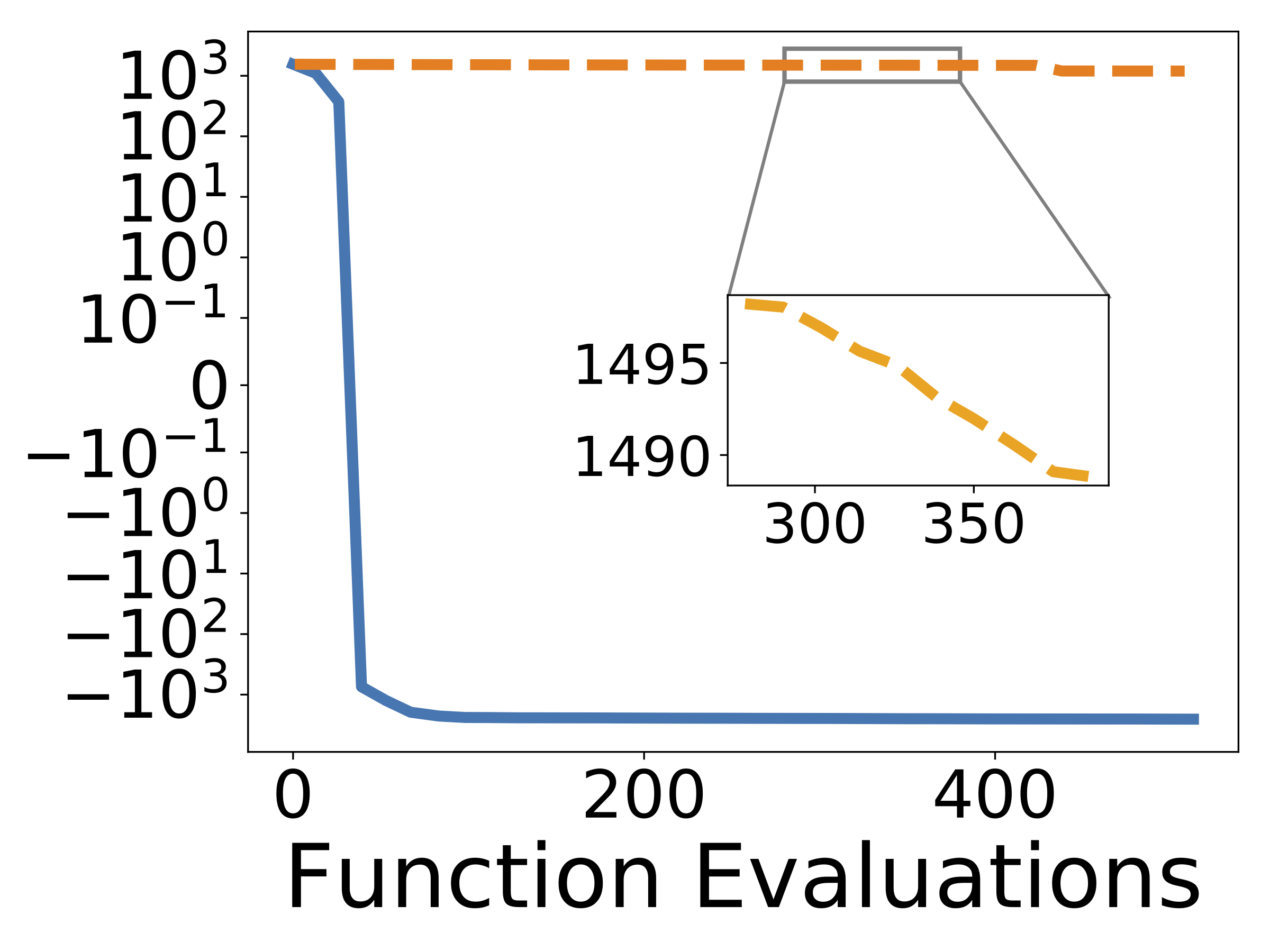}
      \includegraphics[width=.32\linewidth]{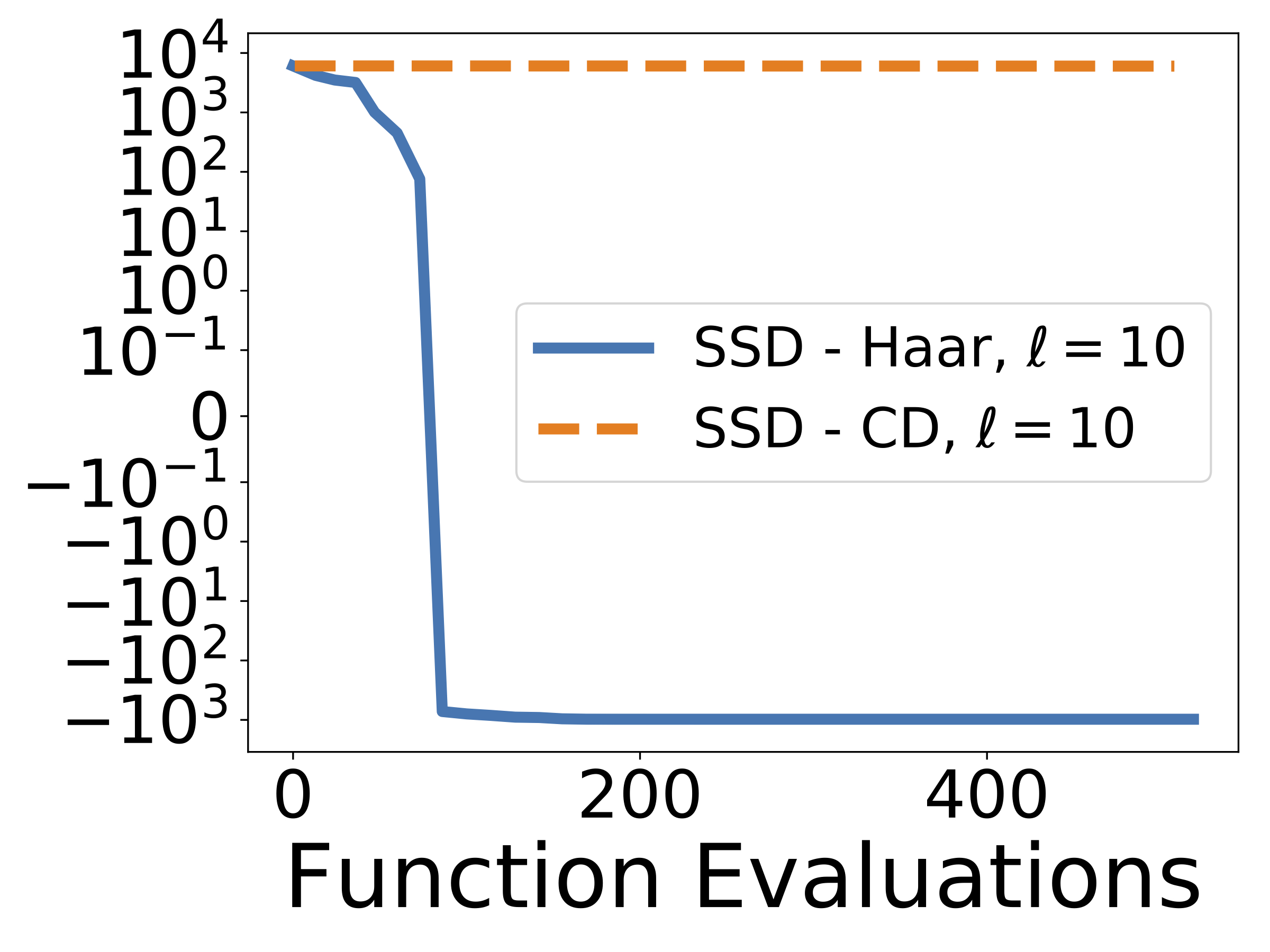}
      \vspace*{-1em}
      \caption{
      Minimizing a function from the family \eqref{eq: NesterovWorst} with $r=d,~\lambda =1$. CD represents randomized block-coordinate descent. Step-size in all cases is chosen by a backtracking linesearch with Armijo conditions. 
      Left: $d=3,~p=50$, total parameters = 153. 
      Center: $d=10,~p=50$, total parameters = 503.
      Right: $d=20,~p=100$, total parameters = 2003. }
      \label{fig: gp}
  \end{figure}
 
The objective function of this problem is non-convex despite the underlying function $T$ being convex. The interpretation of coordinate descent is interesting as each coordinate in parameter space either corresponds to one of the hyperparameters of the kernel, to the noise, or to the placement of one of the inducing points along one dimension. Since $r=d$, the latent function has no low-dimensional structure and movements in any direction in input space correspond to a changing function evaluation. Once again coordinate descent does not scale well with the dimension. This behavior is to be expected: changing the location of particular inducing points along the correct axis has a large improvement on the objective, but if the wrong point is chosen, or the correct point but wrong axis, then little improvement is made (though as we see from the inset, there is slight improvement at each iteration). In contrast, SSD-Haar changes all inducing points in tandem so it descends more rapidly and consistently, particularly in high-dimensional problems. We notice that as before SSD-Haar remains robust to changes in the ambient dimension of the parameter space, though we do see a slight degradation of performance with increased dimension. 
 
\begin{figure}[ht]
      \centering
      \includegraphics[width=.36\linewidth]{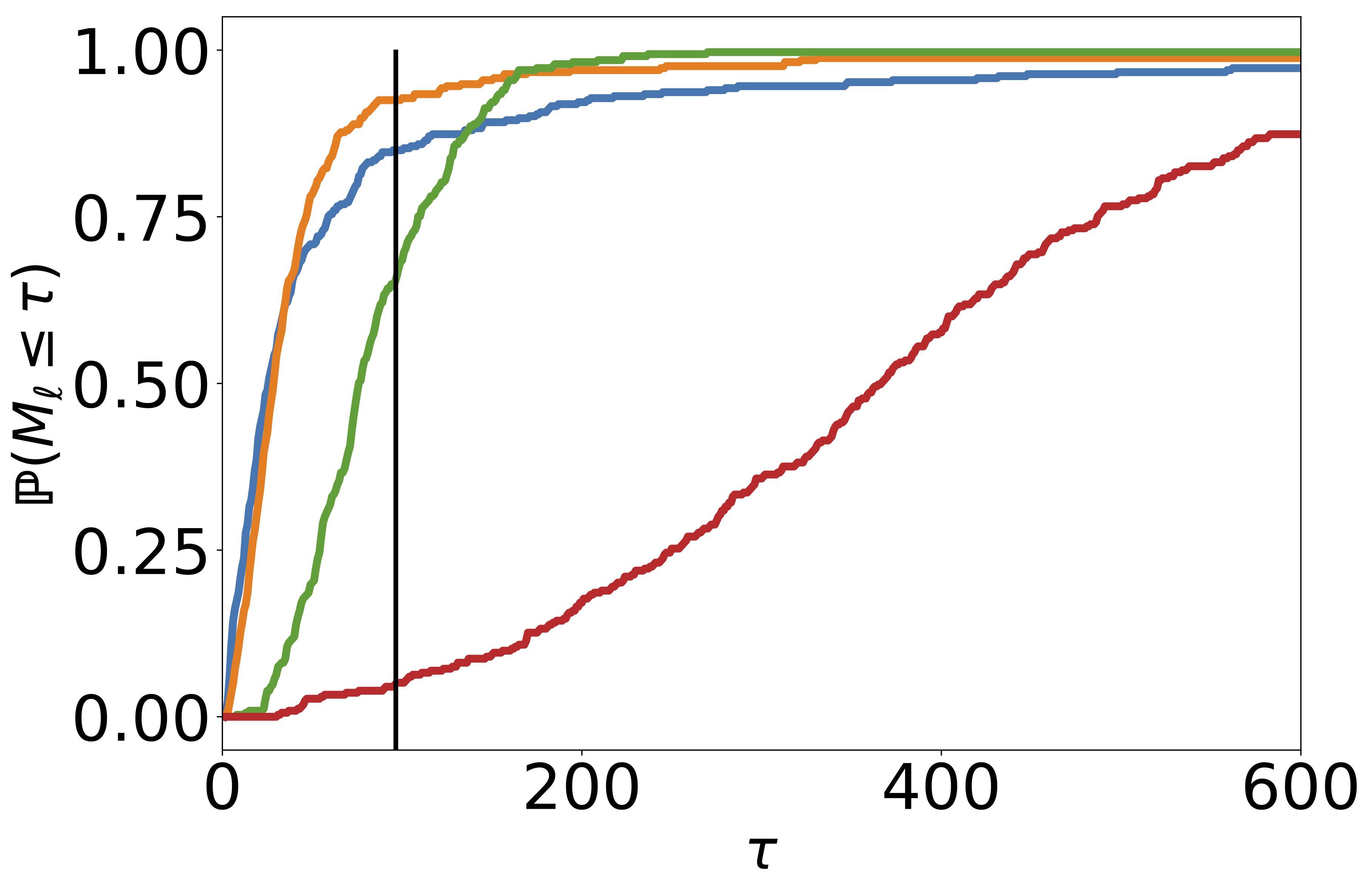}
        \includegraphics[width=.38\linewidth]{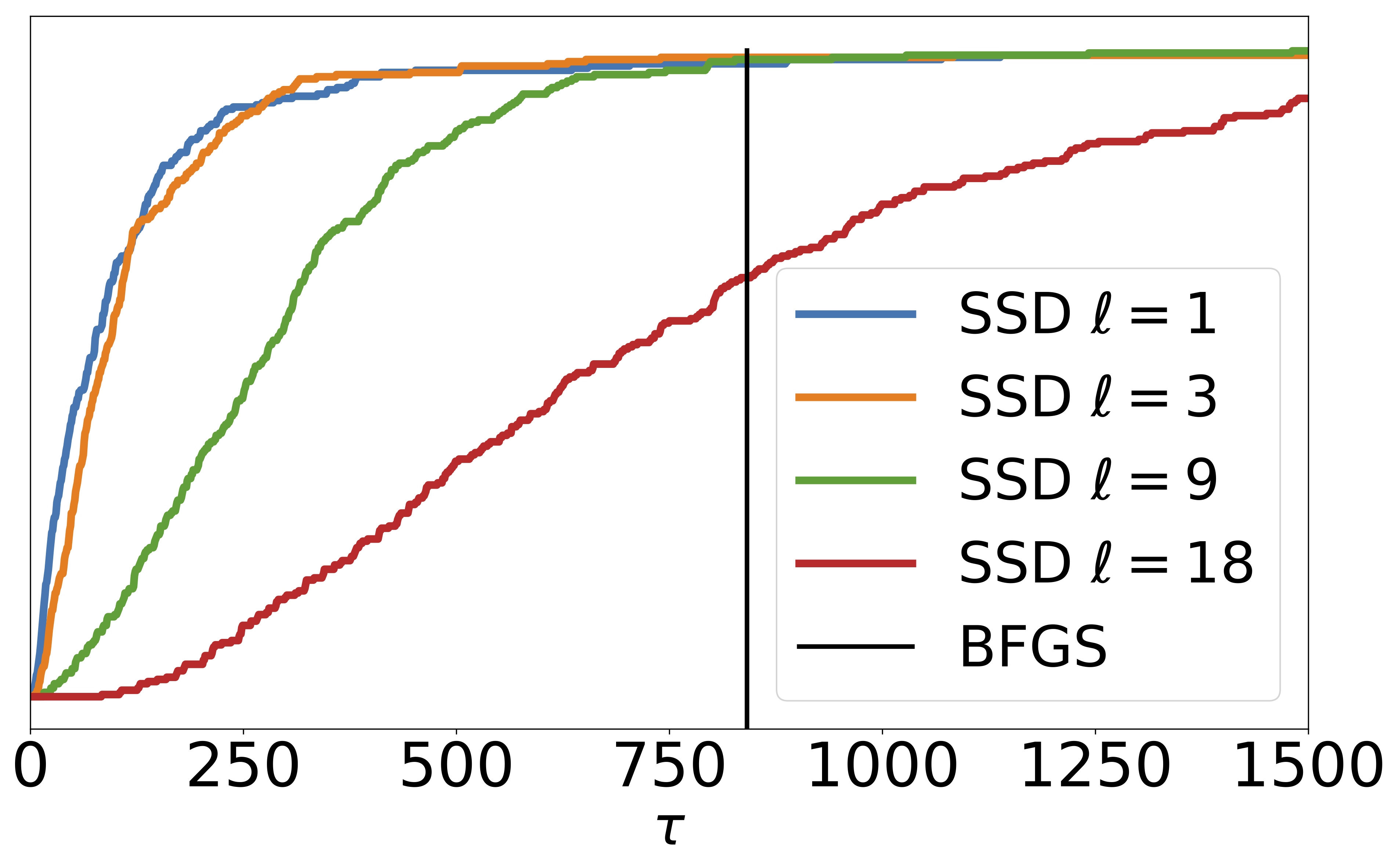}
        \vspace*{-1em}
      \caption{Left: 30-dimensional problem. Right: 60-dimensional problem. $M_{\ell}$ is the number of function evaluations required to attain a cut-off threshold for various values of $\ell$. For a fixed initialization BFGS is non-random, represented by the vertical line. Gradient descent, not pictured, has a vertical line at $\tau=2850$ and $\tau=22828$ for $p=30$ and $p=60$, respectively. $\ell=1$ is equivalent to the method proposed in \cite{nesterov2017random} when $h=0$.}
      \label{fig: SSD-PerformanceProfile}
\end{figure}

  We use performance profiles \cite{dolan2002benchmarking} to determine the effect of varying $\ell$ for different problem sizes and to gauge the variability between runs for a fixed $\ell$. A performance profile is conducted by running each parameterization on a suite of randomized restarts, with termination after some pre-specified tolerance for accuracy has been reached.
  We count the proportion of realizations from each parameterization that achieves the specified tolerance within $\tau$ function evaluations where $\tau=1$ is the fewest function evaluations required in any of the trials, $\tau=2$ is twice as many function evaluations, etc.
  Each parameterization is run 300 times. Results for SSD-Haar are shown in Figure \ref{fig: SSD-PerformanceProfile} for 30- and 60- dimensional objective functions. 

 The cut-off threshold is $95\%$ of the distance between the objective function at the parameter initialization and at the optima, as found by BFGS. Clearly, $\ell=18$ is not a good option in this case. Similarly, $\ell=9$ can be ruled as it underperforms $\ell=1$ and $\ell=3$ approximately 90\% (resp. 99\%) of the time in the 30- (resp. 60-) dimensional problem. The case $\ell=1$ has the best single performance: in the fastest trial it is roughly 100 (resp. 800) times faster than BFGS for the 30- (resp. 60-) dimensional problem, but the variance of the performance for $\ell = 1$ is high, and about 1\% of the time it performs at least 10 times slower than BFGS (not pictured). On the other hand, $\ell=3$ beats BFGS by a similar factor and seems to be insulated from the high variance observed for $\ell=1$. Note also that in 60 dimensions $\ell=3$ is approximately three times faster than BFGS in 90\% of the trials, and about 100 times faster in 40\% of trials. A few trials of $\ell=1$ and $\ell=3$ found their way to a local minima, resulting in the methods not achieving the target threshold.

\subsection{Shape optimization}\label{subsect: experiments-plate}

 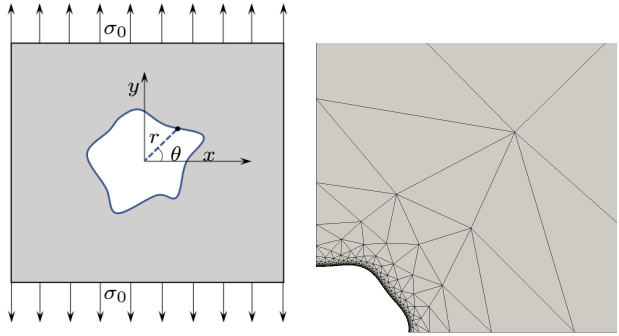
\begin{figure}
      \centering
      \includegraphics[ width=.35\textwidth]{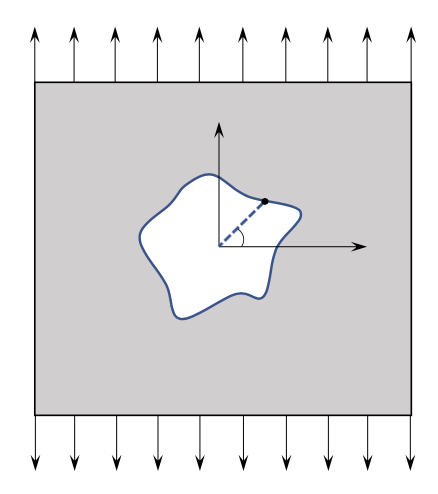}
      \includegraphics[trim=0mm 0mm 0mm 0mm, clip, width=.34\textwidth]{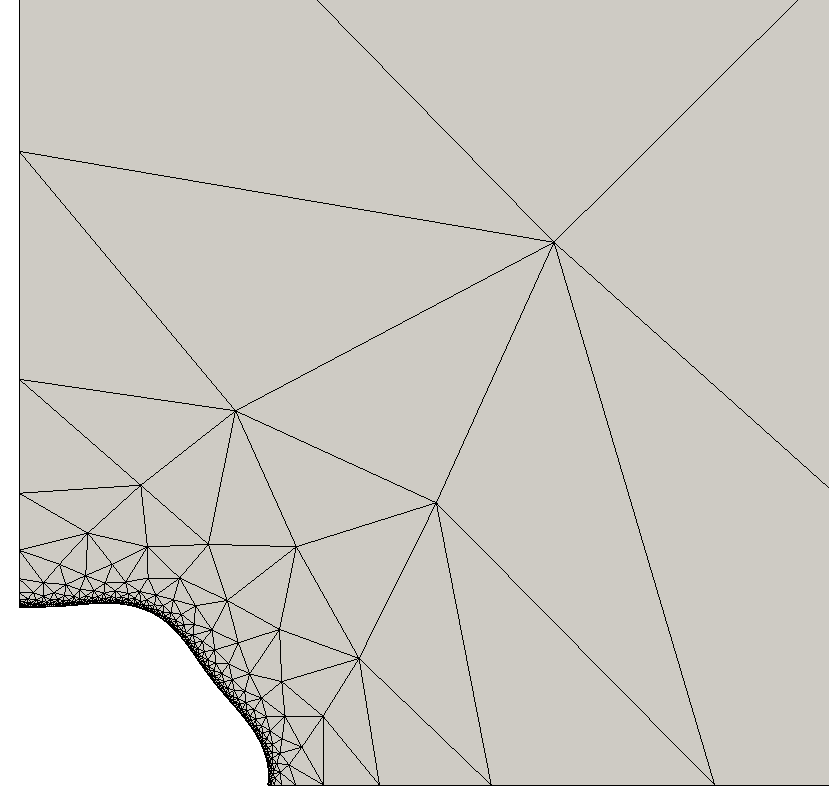}
      \put(-167,67){\footnotesize $\theta$}
      \put(-175,74){\footnotesize $r$}
      \put(-192,115){\footnotesize $\sigma_0$}
      \put(-192,15){\footnotesize $\sigma_0$}
      \put(-155,67){\footnotesize $x$}
      \put(-183,93){\footnotesize $y$}
      \vspace*{-.8em}
      \caption{Left: Schematic of the linear elasticity problem used in the shape optimization example of Section \ref{subsect: experiments-plate}. Right: Conforming finite element mesh used to solve for maximum stress $\sigma_y$ along the $y$ direction. Only a quarter of the plate corresponding to $\theta\in[0,\pi/2]$ is modeled.}
      \label{fig: HoleGeometry}
  \end{figure}
We consider a shape optimization problem involving a linear, elastic structure. Consider a square plate of size $250\times 250$ with a hole, subject to uniform boundary traction $\sigma_0$=1, as illustrated in Fig.~\ref{fig: HoleGeometry}. 
We adopt a discretize-then-optimize approach to solving the PDE-constrained optimization problem. The discretization and optimization steps do not generally commute and an optimize-then-discretize approach  may be preferable for some types of problems \cite[\S 2.9]{gunzburger2003perspectives}, but we do not pursue this question here.

Our goal is to identify a shape of the hole that minimizes the maximum stress $\sigma_y$ along the $y$ direction over a quarter of the plate corresponding to $\theta\in[0,\pi/2]$. To this end, we parameterize the radius of the hole for a given $\theta$ (see Figure \ref{fig: HoleGeometry}) via
 \begin{equation}\label{eqn: parametric-radius}
     r(\theta) = 1 + \delta \medmath{\sum_{i=1}^p} i^{-1/2}\left( \xi_i \sin(i \theta) + \nu_i \cos(i \theta)\right),
 \end{equation}
 where $\delta \in (0,0.5/\sum_{i=1}^p i^{-1/2})$ is a user-defined parameter controlling the potential deviation from an n-gon of radius 1. The parameters that dictate the shape are $\bm{\xi} \in \reals^{p}$ and $\boldsymbol{\nu} \in \reals^{p}$ so that the parameter space is dimension $d = 2p$. Subscripts indicate the index of the vector. We set $\delta = 0.4/\sum_{i=1}^p i^{-1/2}$ so that the minimum possible radius of any particular control point is 0.2 at the initialization. We initialize the entries of $\bm{\xi}$ and $\bm{\nu}$ uniformly at random between -1 and 1. For each instance of $\bm{\xi}$ and $\bm{\nu}$ -- equivalently $r(\theta)$ -- we generate a conforming triangular finite element mesh of the plate that we subsequently use within the FEniCS package \cite{logg2012automated} to solve for the maximum stress $\sigma_y$. A mesh refinement study is performed to ensure the spatial discretization errors are negligible. As we only model a quarter of the plate, we apply symmetry boundary conditions so that $y$ and $x$ displacements along $\theta=0$ and $\theta=\pi/2$ are zero. The Young's modulus and Poisson's ratio of the plate material are set to $E=1000$ and $\nu =0.3$, respectively. A similar problem has been examined in \cite{de2019bi} using a bi-fidelity variant of the popular SVRG algorithm \cite{SVRG}. Due to the  different focus of that work, the investigation of \cite{de2019bi} is conducted in a low-dimensional setting with $d=6$ rather than $d=100$ as in our case.
 
 The parametric radius defined by \eqref{eqn: parametric-radius} enables us to scale the complexity of the problem arbitrarily by increasing the dimension $d$. In effect, if $d$ is large then the problem becomes ill-conditioned since $\xi_p$ and $\nu_p$ each make at most $\delta p^{-1/2}$ additive contribution to the radius. Such ill-posedness suggests that gradient descent ought to perform poorly as it does not account for the curvature of the objective function. Based upon the intrinsic dimensionality results presented in Section \ref{subsect:synthetic-data} we anticipate SSD to outperform gradient descent even though it does not explicitly account for the curvature either. Note that each function evaluation requires a PDE-solve meaning that gradient descent requires $d+1$ PDE-solves per iteration. Though a conforming finite element mesh is used to reduce the computational burden, the cost of so many PDE-solves makes this problem intractable in high-dimensions unless the resolution of the mesh is very low. On the other hand, SSD requires far fewer PDE-solves per iteration provided $\ell \ll d$. 
As mentioned above, the goal is to minimize the maximum stress in the y-direction, $\sigma_y$, over the plate. We make two slight changes to this objective for the sake of the model. First, the stress is obviously minimized if the radius of the hole is zero so we add a term to the objective to penalize deviations from an area of 1 squared unit; even with the regularizer the objective function is still non-convex. Second, the \emph{max} function is not smooth, so it does not fit into the framework of our theory; instead, we minimize the $\ell_p$-norm of the stress with $p=100$, which provides an almost indistinguishable result. 

\begin{figure}[ht]
      \centering
      \includegraphics[width=.32\linewidth]{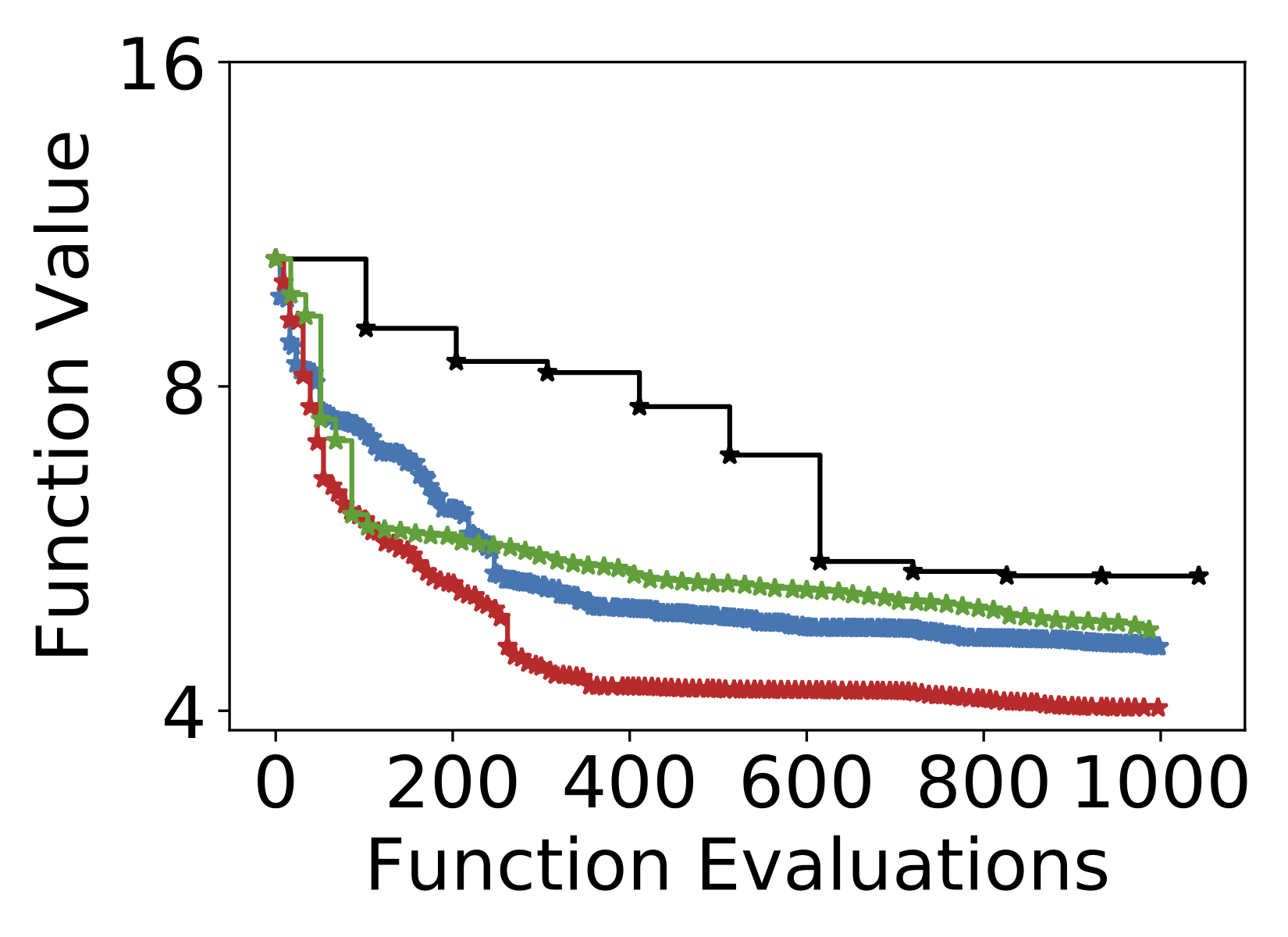}
      \includegraphics[width=.32\linewidth]{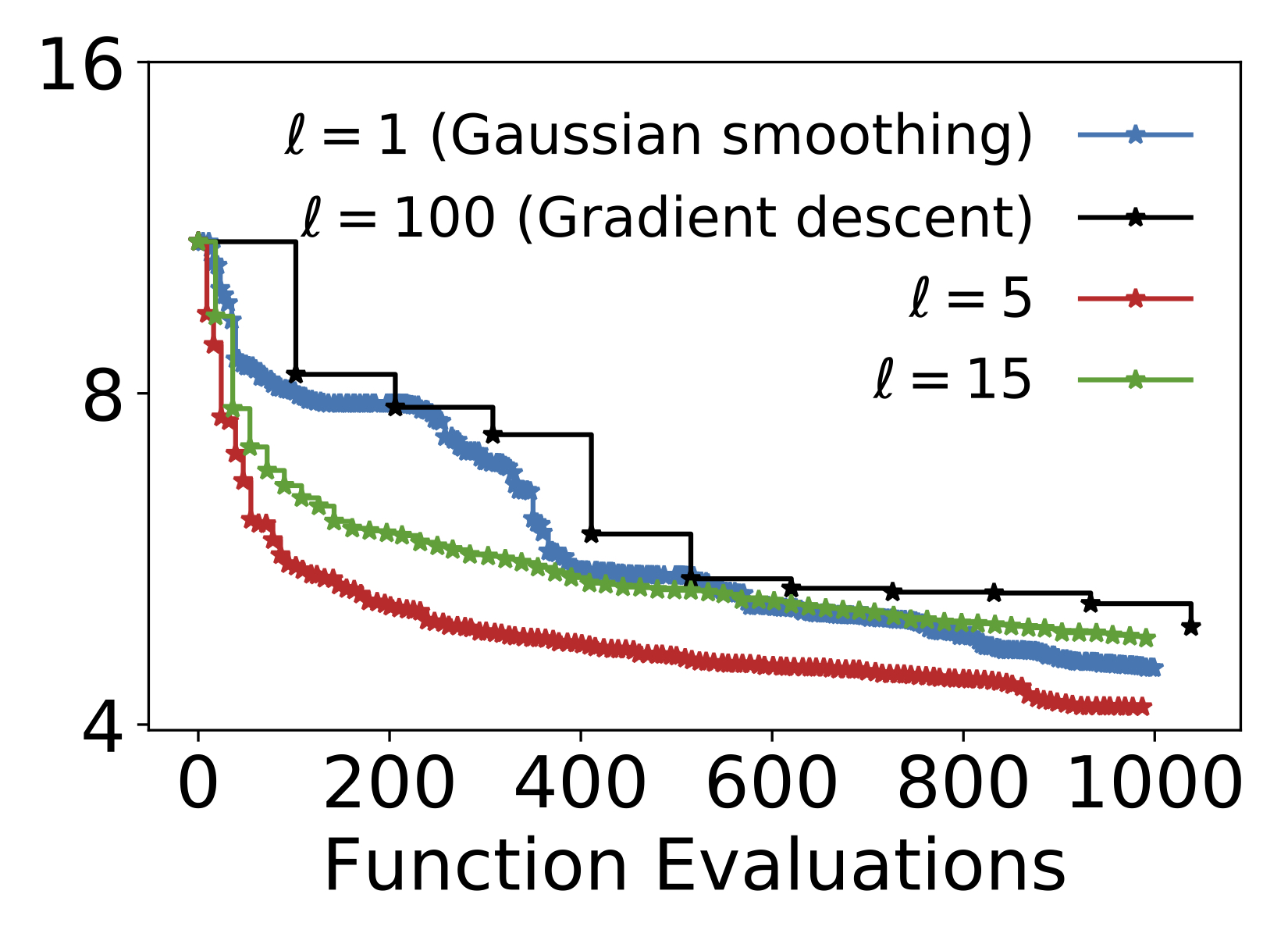}
      \includegraphics[width=.32\linewidth]{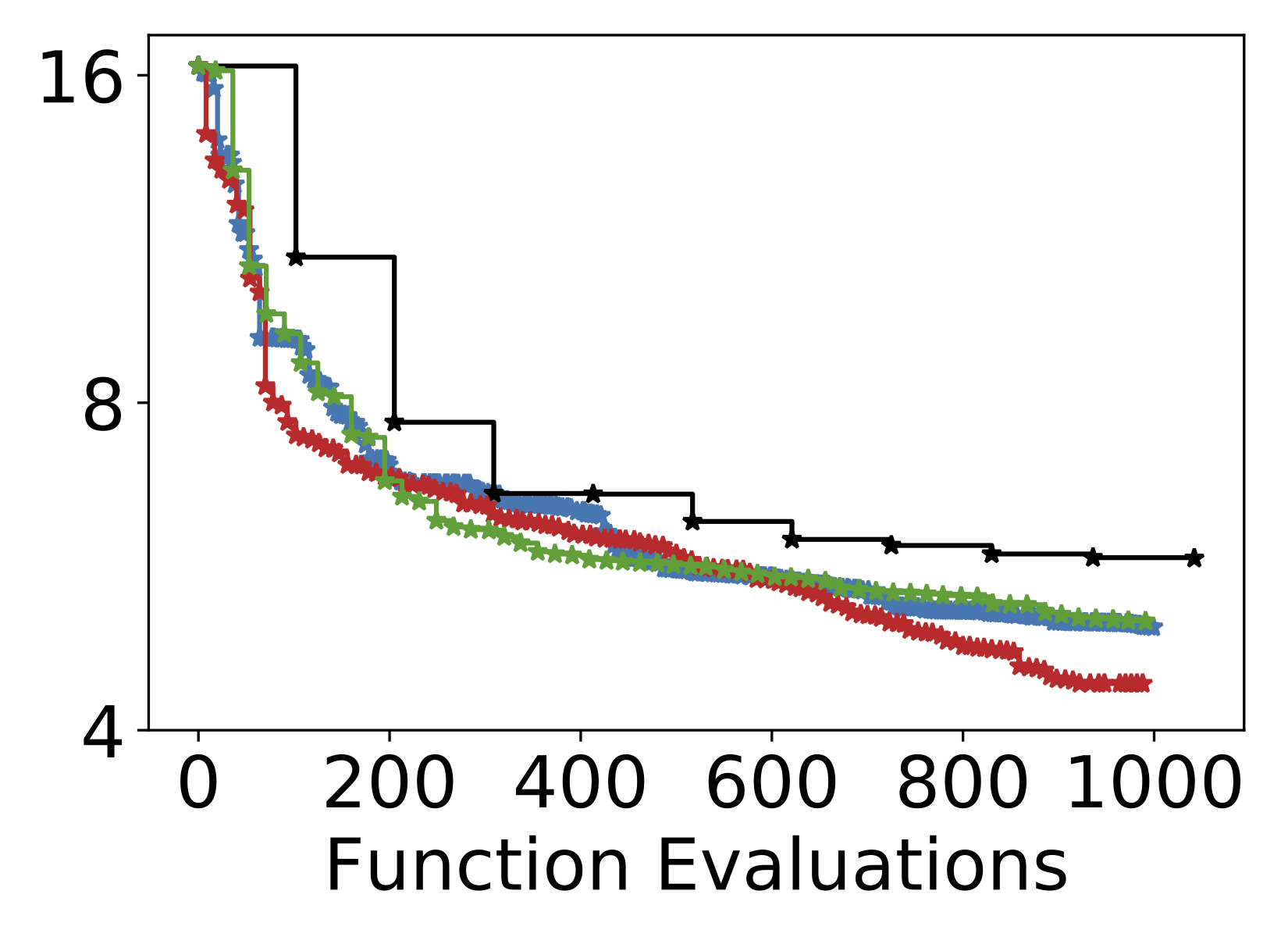}

      \caption{Three runs for optimization of the objective for a hole with shape parameterized by \eqref{eqn: parametric-radius} with $p=50$ (100 dimensions). Each restart represents an initialization of the parameters uniformly at random in $(-1,1)$}.
      \label{fig: SSD-Plate}
  \end{figure}

In Figure \ref{fig: SSD-Plate} we minimize the objective for a hole with shape governed by \eqref{eqn: parametric-radius} for problems with $p=50$ (that is, $d=100$ parameters), using gradient descent and Gaussian smoothing, as well as SSD with $\ell = 5$ and $\ell=15$. In each case, an Armijo backtracking linesearch is used.

In all three randomized restarts finite-difference gradient descent performs poorly relative to the stochastic optimizers. The early iterations are particularly good for the stochastic optimizers. We hypothesize that as the $\ell$-dimensional subspace along which SSD and Gaussian smoothing descends changes with each iteration, parameter space is explored more thoroughly than deterministic methods, making these subspace methods less likely to get funnelled into long, shallow basins; this is intuitively similar to the recent line of research suggesting that noisy perturbation of iterative algorithms helps avoid saddle points \cite{ge2015escaping}. Alternative perspectives hold that subspace methods are cheap on a per-iteration basis so temporarily being caught in a shallow basin is not as expensive in terms of function evaluations. Conversely, a subspace comprised of a single directional derivative (as in Gaussian smoothing) will have a large variance, causing erratic movements through parameter space whenever the gradient is poorly approximated. Figure \ref{fig: SSD-Plate} corroborates the evidence provided in Figure \ref{fig: SSD-PerformanceProfile} that choosing $\ell$ greather than 1 but less than $d$ can be beneficial in terms of rate of convergence.

It is unclear how to choose an optimal $\ell$. Intuition and empirical evidence suggests that a good choice of $\ell$ depends on all of the eigenvalues of $f$, not just on the condition number. In particular, we observe that a rapidly decaying eigenspectrum (as in this problem, and to a larger extent the synthetic data problem described in Section \ref{subsect:synthetic-data}) allows for $\ell$ to be chosen small compared to $d$. In contrast, with a  slow-decaying eigenspectrum choosing $\ell$ small seems to provide relatively less improvement (these experiments are not shown). In none of our experiments does $\ell\ll d$ yield worse results compared to gradient descent when a linesearch is used, suggesting that choosing $\ell \ll d$ may be beneficial with little risk of performing worse. Further analysis must be conducted to verify this assertion. \insertme{An interesting alternative to choosing a fixed $\ell$ is to change $\ell$ as the algorithm progresses in an attempt determine, locally, the appropriate dimension of the subspace used for descent. Further experimentation and analysis would be required to ascertain the benefits of such an adaptive $\ell$.}

\section{Conclusions}\label{sect: Conclusions}
 
We present analysis of an algorithm that generalizes Gaussian smoothing to descend in a randomly chosen subspace and have provided evidence that this generalization is  appropriate for high-dimensional objective functions. We give asymptotic and non-asymptotic results of convergence under a variety of convexity assumptions. We provide tools that are useful beyond the context of this work, such as an interpretation of the Johnson-Lindenstrauss lemma that takes advantage of finite ambient dimension $d$. We demonstrate empirical improvements compared to the \emph{status quo} for several practical problems, and show that the empirical performance can be good even when the assumptions required by the theory are relaxed.

The most obvious extension of this work is a generalization to the case of derivative-free optimization. With directional derivatives unavailable, finite-difference approximations of the derivatives must be employed adding a non-cancelling error at each iteration. Preliminary experiments show that this does not noticeably impede the convergence if $h$, the finite-difference stepsize is sufficiently small.

Thus far, analysis has only been performed for a fixed step-size, but we have shown that an adaptive step-size is required for good practical performance. Recent work in this direction \cite{cartis2018global,berahas2019global} provides promising results that may readily extend to our case. Alternatively, our analysis may be more amenable to trust region methods as in \cite{TrustRegionDFO}

It would be interesting to adapt stochastic optimization algorithms that subsample the observations, as for example in ERM, to the stochastic subspace descent framework. Such sampling would necessitate examination into the effect that noisy function evaluations have on the convergence results. A computationally straightforward extension may allow sketching methods (see e.g. \cite{pilanci2015randomized}) to improve our results with minimal programming overhead, but analysis must be conducted to confirm the theoretical properties of such modifications. 
An adaptive scheme that makes use of observed curvature information could be beneficial for determining the descent directions, an idea that has been discussed at length in the coordinate descent literature \cite{richtarik2014iteration,nesterov2012efficiency}.
Parallelizing our methods to calculate the $\ell$ directional derivatives at each iteration simultaneously is straightforward, but we would like to explore the feasibility of asynchronous parallelization as has been discussed in the coordinate descent case (see, e.g., \cite{peng2016arock}). 
Faster convergence using derivative-free quasi-Newton methods as in \cite{berahas2019derivative} are an obvious extension of this work. 
Finally, recent work on a universal ``catalyst'' scheme \cite{lin2015universal} also applies to our method, allowing for Nesterov-style acceleration without requiring additional knowledge of the Lipschitz constants along any particular direction.

\appendix
\section{Proofs of main results}
\subsection*{Theorem \ref{thm:convergence}} 
	Because $f$ is continuously-differentiable with a $\lambda$-Lipschitz derivative it follows that  
	\begin{equation}\label{Taylor}
	f(\x_{k+1}) \leq  f(\x_k) + \nabla f(\x_k)^\top (\x_{k+1}-\x_k) + \frac{\lambda}{2}\norm{\x_{k+1}-\x_k}^2.
	\end{equation}
Let $f_e(\x) = f(\x)-f_*$ be the error for a particular $\x$.   Then, \eqref{eq: iterations} and \eqref{Taylor} yield:
	\begin{align}\label{delineq}
		f_e(\x_{k+1}) - f_e(\x_k)  &\leq -\alpha_{\lambda} \langle \nabla f(\x_k),~ \mathbf{P}_k\mathbf{P}_k^\top \nabla f(\x_k) \rangle \quad \text{with} \quad \alpha_{\lambda} = \alpha - d\alpha^2\lambda/(2\ell),
	\end{align}
	where we have used the fact that $\mathbf{P}_k\mathbf{P}_k^\top \mathbf{P}_k\mathbf{P}_k^\top = (d/\ell)\mathbf{P}_k\mathbf{P}_k^\top $. Any choice  $0< \alpha < 2\ell/(d\lambda)$ ensures $\alpha_{\lambda} > 0$. With this choice the right hand-side is  non-positive and the errors are non-increasing. Since the error is bounded below by zero the sequence converges almost surely. Furthermore, since the sequence is bounded above by $f_e(\x_0)$, Lebesgue's dominated convergence implies convergence of the sequence in $L^1$. To find the actual limit, define the filtration (i.e., increasing sequence of $\sigma$-algebras) $\mathcal{F}_k = \sigma(\mathbf{P}_1, \ldots, \mathbf{P}_{k-1}), ~ k>1$, and $\mathcal{F}_1 = \{ \emptyset , \Omega\}$. We take conditional expectations of both sides to get
	\begin{equation*}
	\Expectation [f_e(\x_{k+1}) \mid \mathcal{F}_k ] \leq -\alpha_{\lambda} \Expectation [ \langle \nabla f(\x_k), \mathbf{P}_k\mathbf{P}_k^\top  \nabla f(\x_k) \rangle  \mid \mathcal{F}_k ] +f_e(\x_k),
	\end{equation*}
    which leads to
	\begin{equation} \label{intermediate}
	\Expectation \left(f_e(\x_{k+1}) \mid \mathcal{F}_k\right) \leq  -\alpha_{\lambda} \norm{\nabla f(\x_k)}^2 + f_e(\x_{k}),
	\end{equation}
and since $\alpha_{\lambda} > 0$, the PL-inequality yields
	\begin{equation*}
	\Expectation \left(f_e(\x_{k+1}) \mid \mathcal{F}_k\right) \leq -2\gamma  \alpha_{\lambda} f_e(\x_k) + f_e(\x_{k})= \left(1-2\gamma\alpha_{\lambda}\right)f_e(\x_k),
	\end{equation*}
from which we conclude that
	\begin{equation*}
	    	\Expectation f(\x_{k+1})-f_* \leq  \left(1-2\gamma\alpha_{\lambda}\right)^{k+1} \left(f(\x_0) - f_*\right).
	\end{equation*}
    Thus, since $	f_e(\x_k) \convas \X$ for some $\X \in L^1$ and $f_e(\x_k)\overset{L^1}{\longrightarrow} 0$, we have both $f(\x_k) \convas f_*$ and $f(\x_k) \overset{L^1}{\longrightarrow} f_*$. 
\smartqed
\subsection*{Corollary \ref{corr:strong-convexity}$(i)$}
By strong-convexity, the PL-inequality, and Theorem \ref{thm:convergence} we obtain $f(\x_k) \convas f(\x_*)$ and $f(\x_k) - f(\x_*) \geq \frac{\gamma}{2} \norm{\x_*-\x_k}$. Since the left-hand side converges a.s.\ to zero and $\gamma > 0$, we have $\x_k \convas \x_*$.
\subsection*{Corollary \ref{corr:strong-convexity}$(ii)$}
 Rearranging the terms in equation \eqref{intermediate} we have \\$-\alpha_{\lambda} ^{-1}\Expectation \left(f_e(\x_k) - f_e(\x_{k+1})\mid \mathcal{F}_k \right) \geq \norm{\nabla f(\x_k)}^2$. Combining this with Lipschitz continuity yields
	$
	2\gamma f_e(\x_k) \leq \norm{\nabla f(\x_k)}^2 \leq -\alpha_{\lambda}^{-1} \Expectation\left(f_e(\x_{k}) - f_e(\x_{k+1}) \mid \mathcal{F}_k\right). $
	That is,
	\begin{equation}\label{nearly}
	\Expectation \left(f_e(\x_{k+1}) \mid \mathcal{F}_k\right) \leq \left(1-2\gamma\alpha_{\lambda}\right)f_e(\x_k).
	\end{equation}
	Choosing $\alpha_{\lambda}=\ell/(d\lambda)$ results in $\Expectation f_e(\x_{k+1}) \leq \left(1-\ell\gamma / (d \lambda) \right)^{k+1}f_e(\x_0)$
\smartqed
	\subsection*{Theorem \ref{thm: convergence-convex}}
	We follow the proof of Theorem \ref{thm:convergence} until \eqref{intermediate}, then we rearrange terms to obtain,
	\begin{align}\label{eq:simplified-Lipschitz}
	\Expectation \left(f(\x_{k+1}) \mid \mathcal{F}_k \right) \leq f(\x_k) - \alpha_{\lambda} \norm{\nabla f(\x_k)}^2,
	\end{align}
	and then by convexity and the Cauch-Schwarz inequality, $\norm{\nabla f(\x_k)} \geq f_e(\x_k)/R$. Plugging this into equation \eqref{eq:simplified-Lipschitz} and letting $\alpha = \ell/(d\lambda)$ results in 
	\begin{equation}\label{eq: required}
 \Expectation [f_e(\x_{k+1}) \mid \mathcal{F}_k]- f_e(\x_k)  \leq  -\alpha  f_e(\x_k)^2/2 R^2,
	\end{equation}
	and one more expectation yields
	\begin{align*}
	    \Expectation[f_e(\x_{k+1}) - f_e(\x_k)] \leq -\alpha \Expectation f_e(\x_k)^2/2R^2
	    &\leq -\alpha \left(\Expectation f_e(\x_k)\right)^2/2R^2 \\
	    \leq -\alpha \Expectation f_e(\x_k)\cdot \Expectation f_e(\x_{k+1})/(2R^2)
	\end{align*}
	since $\alpha\ge 0$ and $\Expectation f_e(\x_{k+1}) \le \Expectation f_e(\x_k)$. Dividing by $\Expectation f_e(\x_k)\cdot \Expectation f_e(\x_{k+1})$ gives
	%
	\begin{equation}\label{eq: before recursion}
	\frac{1}{\Expectation f_e(\x_{k+1})} \geq 	\frac{1}{\Expectation f_e(\x_{k})} + \frac{\alpha}{2 R^2}.
	\end{equation}
	Applying \eqref{eq: before recursion} recursively, and replacing $\alpha$ with $\ell/(d\lambda)$ we obtain
	$\Expectation  f_e(\x_{k+1}) \leq 2d\lambda R^2/k \ell$. 
\smartqed
\subsection*{Theorem \ref{thm: convergence-nonconvex}}
Beginning from \eqref{intermediate} we set $\alpha_{\lambda} = \ell/(d\lambda)$ and rearrange terms to get
\begin{equation*}
\ell/(2d\lambda) \norm{\nabla f(\x_k)}^2 \leq f(\x_k) - \Expectation (f(\x_{k+1})\mid \mathcal{F}_k),
\end{equation*}
which leads to
\begin{equation*}
    \ell/(2d\lambda) \sum_{i=0}^k \Expectation \norm{\nabla f(\x_k)}^2 \leq \sum_{i=0}^k \Expectation (f(\x_i) - f(\x_{i+1})) = f(\x_0) - \Expectation f(\x_{k+1}) \leq f(\x_0)-f_*. 
\end{equation*}
Recognizing that a sum of $k+1$ values is bounded below by $k+1$ replicates of its minimum yields
\begin{equation*}
    (k+1) \min_{i\in\{0,\dots,k\}} \Expectation \norm{\nabla f(\x_k)}^2 \leq \frac{2d\lambda (f(\x_0)-f_*)}{\ell}.
\end{equation*}
Divide both sides by $k+1$ to get the result. Now, define some tolerance $\epsilon$ such that 
\begin{equation*}
    \frac{2d\lambda (f(\x_0)-f_*)}{(k+1)\ell} \leq \epsilon.
\end{equation*}
Then, 
\begin{equation*}
    k \geq  \frac{2d\lambda (f(\x_0)-f_*)}{\epsilon\ell}-1.
\end{equation*}
That is, $k = \mathcal{O}(d/\ell\epsilon)$ iterations are sufficient to achieve $\Expectation\norm{\nabla f(\x_k)} \leq \epsilon$.

\subsection*{Lemma \ref{Lemma: DavidJL}}
Let  $\mathbf{H} \in \reals^{d \times d}$ be a Haar-distributed random matrix, $\v \in \reals^d$ an arbitrary fixed vector, and $\mathbf{u} \sim \mathcal{N}(\mathbf{0},\I_d)$. 
Then $\mathbf{H}^\top \v/\norm{\v}$  and $\mathbf{u}/\norm{\mathbf{u}}$ are both distributed uniformly on the $d$-dimensional sphere.
Let $\I_{\ell \times d} \in \reals^{\ell \times d}$ represent a mapping onto the first $\ell$ coordinates. Then, 
\begin{equation*}
\norm{\I_{\ell \times d}\mathbf{u}}^2 =\left(u_1^2 + \ldots + u_\ell^2\right) \sim \chi^2(\ell),
\end{equation*} 
and 
\begin{align}\label{eqn: chisq}
\nonumber \norm{\mathbf{u}}^2 = u_1^2 + \ldots + u_\ell^2 + u_{\ell+1}^2 + \ldots + u_d^2 &\sim \chi^2(d).
\end{align}
For independent random variables $X \sim \chi^2(\alpha)$ and $Y \sim \chi^2(\beta)$, $Z = X/(X+Y) \sim \mathcal{B}\mathrm{eta}(\alpha/2, \beta/2)$. Thus,
\begin{equation*}
\frac{\norm{\I_{\ell \times d}\mathbf{H}^\top \v}^2}{\norm{\v}^2}=\norm{\I_{\ell \times d}\mathbf{H}^\top \frac{\v}{\norm{\v}}}^2\overset{d}{=}\norm{\I_{\ell \times d} \frac{\mathbf{u}}{\norm{\mathbf{u}}}}^2 =\frac{\norm{\I_{\ell \times d} \mathbf{u}}^2}{\norm{\mathbf{u}}^2}\sim \mathcal{B}\mathrm{eta}(\ell/2, (d-\ell)/2).
\end{equation*}
By construction, $\P_k\overset{d}{=}\sqrt{d/\ell}\,\I_{\ell \times d}\mathbf{H}$, so
 \begin{align*}
     \mathbb{P}\left(\norm{\P_k^\top\v}^2 \leq (1-\epsilon)\norm{\v}^2  \right) &= \mathbb{P}\left( \norm{\I_{\ell \times d}\mathbf{H}^\top \frac{\v}{\norm{\v}}}^2 \leq  \frac{\ell}{d} (1-\epsilon) \right).
 \end{align*} 
 The Beta CDF is calculated by evaluating the regularized incomplete Beta function. That is, if $X \sim \mathcal{B}eta(\alpha,\beta)$ then $F_X(p) = I_p(\alpha,\beta)$. Thus, the probability
 \begin{equation*}
     \mathbb{P}\left( \norm{\I_{\ell \times d}\mathbf{H}^\top \frac{\v}{\norm{\v}}}^2 \geq  \frac{\ell}{d} (1-\epsilon) \right) =1-I_{(1-\epsilon)\ell/d}(\ell/2, (d-\ell)/2)
 \end{equation*}
provides a probability of a successful embedding.
\insertme{
\subsection*{Remark \ref{remark: independence}}
We show that for any $ k_1 < \cdots < k_m$ the sets     $A_{k_1},\ldots, A_{k_m}$
are mutually independent. Let $\mathds{1}_A$ denote the indicator function of a set $A$. Define the filtration (i.e., increasing sequence of $\sigma$-algebras) $\mathcal{F}_k = \sigma(\mathbf{P}_1, \ldots, \mathbf{P}_{k-1}), ~ k>1$, and $\mathcal{F}_1 = \{ \emptyset , \Omega\}$ and note that $A_k$ is $\mathcal{F}_{k-1}-$measurable. Then by the chain rule of probability and the fact that the $(\P_k)$ are iid,
\begin{align*}
    \Prob(A_{k_1} \cap \cdots \cap A_{k_m}) &= \Expectation [\mathds{1}_{A_{k_1}}\cdots\mathds{1}_{A_{k_{m-1}}}]\Prob(A_{k_m} \mid \mathcal{F}_{k_{m-1}}) \\
    &= \Expectation [\mathds{1}_{A_{k_1}}\cdots\mathds{1}_{A_{k_{m-2}}}]\Prob(A_{k_{m-1}} \mid \mathcal{F}_{k_{m-2}})  \Prob(A_{k_{m}} \mid \mathcal{F}_{k_{m-1}}) \\
    &\,\vdots \\
    &= \Prob(A_{k_1} \mid \mathcal{F}_{k_0}) \cdots \Prob(A_{k_{m}} \mid \mathcal{F}_{k_{m-1}}) \\
    &= \Prob(A_{k_1}) \cdots \Prob(A_{k_m}).
\end{align*}
}
\subsection*{Theorem \ref{thm: strong-convexity}}
Beginning from \eqref{delineq} we choose an $\ell$ determined by Lemma \ref{Lemma: DavidJL} such that with probability $\delta$,
	\begin{equation}\label{eqn: intermediate}
	    f_e(\x_{k})  \leq  f_e(\x_{k-1}) - (1-\epsilon)\alpha_{\lambda} \norm{\nabla f(\x_{k-1})}^2.
	\end{equation}
    By (A3') the function is $\gamma$-strongly-convex, so, 
	\begin{equation}\label{eq: after-strong-convexity}
  f_e(\x_{k}) \leq  \left(1-(1-\epsilon)\frac{\ell \gamma}{d \lambda}\right)f_e(\x_{k-1})\quad \text{with probability }\delta.
	\end{equation}
    Define the Bernoulli random variable $W_k \sim Bern(\delta)$ such that $W_k=1$, occurring with probability $\delta$, constitutes a successful embedding on the $k^\text{th}$ iteration. We can re-write \eqref{eq: after-strong-convexity} as
    \begin{equation}\label{eqn: using indicators}
        f_e(\x_{k}) \leq \left(1 - W_k(1-\omega)\right)f_e(\x_{k-1}),
    \end{equation}
    where $\omega = 1-(1-\epsilon)\ell\gamma/(d\lambda)$. If the embedding is a failure, we use the trivial bound $\norm{\P_k^\top \nabla f(\x_k)}^2\geq 0$. Consider a random variable $U_k =1- W_k(1-\omega)$, then
    \eqref{eqn: using indicators} is
    \begin{equation*}
        f_e(\x_{k}) \leq (U_1\cdots U_k)f_e(\x_0).
    \end{equation*}
    Note that $\log(U_k) = Y_k \log \omega$ for  $Y_k \sim \text{Bernoulli}(\delta)$. 
    Let $B \sim \text{Bin}(k,\delta)$, then, for $t' \in(0,\,k\delta] $
\begin{equation}\label{eq: not-perp}
    \Prob(U_1\cdots U_k \geq \omega^{k\delta-t'}) = \Prob(B \log \omega \geq (k\delta-t')\log \omega) = \Prob(B \leq k\delta-t').
\end{equation}
Thus, for $t' \in (0,k\delta]$ we obtain a probabilistic lower bound on the improvement using Remark \ref{Lemma: Optimal binomial},
\begin{equation*}
\Prob( U_1 \cdots U_k \geq \omega^{k \delta  - t'}) \leq \exp(-t'^2/2\sigma_k^2),
\end{equation*}
where $\sigma^2_k = k(1-2\delta) / (2\log((1-\delta)/\delta))$. Now,
\begin{align*}
     \Prob \left( f_e(\x_{k}) \geq \omega^{k \delta  - t'}  \right) &\leq \Prob \left((U_1 \cdots U_k) f_e(\x_0) \geq \omega^{k \delta  - t'}  \right) \\ 
     &= \Prob((U_1 \cdots U_k) \geq \omega^{k \delta  - t'} /f_e(\x_0))
\end{align*}
which implies that for $t' \in(0,\,k\delta]$,
\begin{equation*}
\Prob\left(f_e(\x_{k}) \geq \left(1-(1-\epsilon)\frac{\ell \gamma}{d \lambda}\right)^{k \delta  - t'}f_e(\x_0) \right) \leq \exp(-t'^2/2\sigma_k^2).
\end{equation*}
Define $t = (t'/k) \in (0,\delta]$ and the result follows.

    \smartqed

\begin{acknowledgements}
We thank Lior Horesh for his suggestions regarding derivative-free optimization, Gregory Fasshauer for fruitful discussions on sparse Gaussian processes, and Neil Longfellow and Osman Malik for their insights into connections with the Johnson-Lindenstrauss lemma. We also thank Subhayan De for providing the FEniCS code used in Section \ref{subsect: experiments-plate}. AD acknowledges funding by the US Department of Energy’s Office of Science Advanced Scientific Computing Research, Award DE-SC0006402 and National Science Foundation Grant CMMI-145460. LT acknowledges funding by National Science Foundation grant DMS-1723005. SB acknowledges funding by National Science Foundation grant DMS-1819251.
\end{acknowledgements}

\bibliographystyle{siam}
\bibliography{RandomizedBib,thesisBecker}

\end{document}